# A Fast Algorithm for Convolutional Structured Low-rank Matrix Recovery

Greg Ongie*, *Member, IEEE*, Mathews Jacob, *Senior Member, IEEE*

**Abstract**

Fourier domain structured low-rank matrix priors are emerging as powerful alternatives to traditional image recovery methods such as total variation and wavelet regularization. These priors specify that a convolutional structured matrix, i.e., Toeplitz, Hankel, or their multi-level generalizations, built from Fourier data of the image should be low-rank. The main challenge in applying these schemes to large-scale problems is the computational complexity and memory demand resulting from lifting the image data to a large scale matrix. We introduce a fast and memory efficient approach called the Generic Iterative Reweighted Annihilation Filter (GIRAF) algorithm that exploits the convolutional structure of the lifted matrix to work in the original un-lifted domain, thus considerably reducing the complexity. Our experiments on the recovery of images from undersampled Fourier measurements show that the resulting algorithm is considerably faster than previously proposed algorithms, and can accommodate much larger problem sizes than previously studied.

## I. Introduction

Recovering image data from limited and noisy measurements is a central problem in many biomedical imaging problems, including MRI, microscopy, and tomography. The conventional approach to solve these ill-posed problems is to regularize the recovery using smoothness or sparsity priors in discrete image domain, as in compressed sensing-based MRI reconstruction [1]. In contrast to this approach, several researchers have recently proposed a new class of recovery algorithms that impose a low-rank penalty on a convolutional structured matrix constructed from Fourier data of the image, i.e., a matrix that is Toeplitz, Hankel, or their multi-level generalizations [2]–[8]. Structured low-rank matrix priors rely on linear dependencies among the Fourier coefficients, which arise from variety of assumptions including continuous domain analogs of sparsity [3], [6], [8], [9], correlations in the locations of the sparse coefficients in spatial domain [8], [9], multi-channel sampling [2], [10], [11], or smoothly varying complex phase in spatial domain [4]. For example, we have shown that the Fourier coefficients of continuous domain piecewise constant images whose discontinuities are localized to zero level-set of a bandlimited function satisfy an annihilation relationship [9]; their recovery from undersampled Fourier measurements translates to a convolutional structured low-rank matrix completion problem [8]. This model, which exploits the smooth structure of the discontinuities along with sparsity, can lead to significant improvement in reconstruction quality over traditional methods such as total variation (TV) minimization; see [8], [9] for more details. Similarly, the ALOHA framework [6] reformulates the recovery of a transform sparse signal as a structured low-rank matrix recovery problem, e.g., images that are sparse under the undecimated Harr wavelet transform. These convolutional structured low-rank matrix penalties are emerging as powerful alternatives to conventional discrete spatial domain penalties because of their ability to exploit continuous domain analogs of sparsity. These methods are often called "off-the-grid", in the sense that they do not

G. Ongie is with the Department of Electrical Engineering and Computer Science, University of Michigan, Ann Arbor, MI, 48109, and M. Jacob is with the Department of Electrical and Computer Engineering, University of Iowa, Iowa City, IA, 52245 USA (e-mail: gongie@umich.edu; mjacob@uiowa.edu)

This work is supported by grants NIH 1R01EB019961-01A1, NSF CCF-1116067, ONR N00014-13-1-0202, and ACS RSG-11-267-01-CCE.

require discrete gridding of the signal in spatial domain. This is because these methods model linear dependencies among the Fourier coefficients of the image induced by low-complexity continuous domain properties of the image (detailed below). In particular, the linear dependencies can be expressed as a nulling of the Fourier coefficients of the image by convolution with a finite impulse response filter—a so-called "annihilation relationship"—that translates directly to the low-rank property of a Toeplitz-like matrix built from the Fourier coefficients. These methods then pose image recovery either partially or entirely in Fourier domain as a structured low-rank matrix recovery problem, where the structure of the matrix encodes the particular modeling assumptions. The LORAKS framework [4] capitalized on the discrete domain sparsity and smooth phase of the continuous domain image using structured low-rank matrix priors, which offers improved reconstructions over conventional $\ell_1$ recovery. Inspired by auto-calibration techniques in parallel MRI [12]–[15], the extension to recovery of parallel MRI from undersampled measurements is formulated as a structured low-rank matrix recovery problem in [2], [10]. Similar approaches have also been found very effective in auto-calibrated multishot MRI [16] and correction of echo-planar MRI data [17]. The theoretical performance of structured low-rank matrix completion methods has been studied in [3], [18], [19], showing improved statistical performance over standard discrete spatial domain recovery.

Despite improvements in empirical and statistical recovery performance over traditional methods as seen from [2]–[4], [6], [8], [16], these structured low-rank matrix recovery schemes are associated with a dramatic increase in computational complexity and memory demand; this restricts their direct extension to multi-dimensional imaging applications, such as dynamic MRI reconstruction. In particular, these algorithms involve the recovery of a large-scale Toeplitz-like matrix whose combined dimensions are several orders of magnitude larger than those of the image. For example, the dimension of the structured matrix is roughly $10^6 \times 2000$ for a realistic scale 2-D MRI reconstruction problem (see Table I). The analogous matrix lifting is several orders of magnitude larger in dynamic MRI applications, making it impossible to store. Moreover, typical algorithms require a singular value decomposition (SVD) of the dense large-scale matrix at each iteration, which is computationally prohibitive for large-scale problems. Several strategies have been introduced to minimize or avoid SVD's in these algorithms. For example, the algorithms derived in [4], [6] replace the full SVD's with more efficient truncated SVD's or matrix inversions by assuming the matrix to be recovered is low-rank, or well approximated as such. However, even with these low-rank modifications, the algorithms in [4], [6] still have considerable memory demand, since they require storing a variable having dimensions of the large-scale matrix.

In this paper we introduce a novel, fast algorithm for a class of convolutional structured low-rank matrix recovery problems arising in MRI reconstruction and other imaging contexts. What distinguishes the algorithm from other approaches is its direct exploitation of the convolutional structure of the matrix—none of the current algorithms fully exploit this structure. This enables us to evaluate the intermediate steps of the algorithm efficiently using fast Fourier transforms in the original problem domain, resulting in an algorithm with significant reductions in memory demand and computational complexity. Our approach does not require storing or performing computations on the large-scale structured matrix, nor do we need to make overly strict low-rank assumptions about the solution. The proposed approach is based on the iterative reweighted least-squares (IRLS) algorithm for low-rank matrix completion [20]–[22], which we adapt to the convolutional structured matrix setting. This algorithm minimizes the nuclear norm or the non-convex Schatten-$p$ quasi-norm of the structured matrix subject to data constraints. However, as we show, the direct extension of the IRLS algorithm to the large-scale setting does not yield a fast algorithm. We additionally propose a systematic approximation of convolutional structured matrices that radically simplifies the subproblems of the IRLS algorithm, while keeping the low-rank property of the matrix intact. The combination of these two ingredients, the IRLS algorithm and the proposed approximation of the convolutional structured matrices, we call the generalized iterative reweighted annihilating filter (GIRAF) algorithm. The name reflects the fact that the algorithm can be interpreted as alternating between (1) the robust estimation of an annihilating filter for the data, and (2) solving for the data best annihilated by the filter in a least-squares sense.

The GIRAF algorithm can also be viewed as a generalization of the two-step "auto-calibrated" recovery schemes used in many structured low-rank models for signal recovery. These schemes first solve for one or many annihilating filters using a non-iterative approach from a uniformly sampled *calibration region* in Fourier domain, which is then used to interpolate missing values in the remaining Fourier data by solving a least-squares linear prediction problem. For example, this approach is commonly used in GRAPPA and related techniques in parallel MRI reconstruction [12]–[15], as well as in AC-LORAKS [23], which was proposed to speed up MRI reconstruction using LORAKS modeling [4], [10] when calibration data is present. Similarly, a two-step approach was adopted for super-resolution MRI recovery [9], [24]. The main benefit of the GIRAF scheme is its ability to work in variety of sampling settings, including the "calibrationless" setting, where there in no fully sampled calibration region. In contrast to these works, the GIRAF algorithm allows us to obtain an annihilating filter from *non-uniform* samples in Fourier domain via an iterative procedure, which also recovers the missing Fourier data by solving a similar weighted least-squares problem at each iteration.

Finally, we note that several existing approaches for inverse problems in imaging utilize non-convex approaches, similar to the present work. For example, an augmented Lagrangian iterative shrinkage algorithm for non-convex Schatten $p$-norms has been introduced for the denoising of images in [25]. However, since these schemes apply low-rank denoising on groups of non-local patches, the matrices in their setting are considerably smaller than in our context and does not possess the convolutional structure that we consider in this work. Similarly, the non-convex penalties studied in [26], [27] are discrete formulations using non-convex vector norms in the image domain; these methods cannot exploit the continuous domain properties that are modeled using structured low-rank penalties, and hence will have limited applicability in the context we consider.

## II. Signal reconstruction by convolutional structured low-rank matrix recovery

### A. *Convolutional structured low-rank models in imaging*

Structured low-rank matrix approximation (SLRA) models are widely used in many branches of signal processing [28], [29]. A SLRA model assumes some property of the signal data is equivalent to the low-rank property of a matrix constructed from the data. This paper is motivated by recent convolutional SLRA models in MRI reconstruction [2], [4]–[8], [24], and related inverse problems in imaging [6], [30], [31]. In this setting, various spatial domain properties of the image (e.g., limited support, smooth phase, piecewise constant, etc.) translate into the low-rank property of a convolutional structured matrix (e.g., Toeplitz, Hankel, and their multi-level generalizations) constructed from the Fourier coefficients of the image. Recovery of the image from undersampled or corrupted measurements is then posed in Fourier domain as a convolutional structured low-rank matrix recovery problem.

As a motivating example for the SLRA approach in imaging, consider the class of signals consisting of a sparse linear combination of Dirac impulses in 1-D:

$$\rho(x) = \sum_{i=1}^{r} c_i\, \delta(x - x_i), \quad \text{for all} \quad x \in [0,1] \tag{1}$$

This signal model is "off-the-grid" in the sense that the impulse locations $\{x_i\}_{i=1}^r$ can be arbitrary points in $[0,1]$ and are not required to lie on a discrete grid. This example and the analysis that follows is closely related classical methods in line spectral estimation, including Prony's method [32] and its robust refinements: MUSIC [33], ESPRIT [34], matrix pencil [35], and others. See [36] for a comprehensive overview. These methods utilize the fact that the sparsity of the signal implies a Toeplitz matrix built from the Fourier coefficients of $\rho$ is rank deficient. To see why this is, let $\mu(x)$ be a periodic bandlimited function on $[0,1]$ having $r$ zeros at the Dirac locations $\{x_i\}_{i=1}^r$. It is easily seen that

$$\rho(x)\mu(x) = 0 \quad \text{for all} \quad x \in [0,1] \tag{2}$$

where the equality is understood in the sense of distributions or generalized functions (see, e.g., [37]). This multiplication annihilation relationship in spatial domain translates to a convolution annihilation relationship in Fourier domain:

$$(\hat{\rho} * \hat{\mu})[k] = 0, \quad \text{for all} \quad k \in \mathbb{Z} \tag{3}$$

where $\hat{\rho}$ and $\hat{\mu}$ denote the Fourier coefficients of $\rho$ and $\mu$. For a finite collection of low-pass Fourier coefficients $\hat{\rho}[k], |k| \leq K$, where $K$ is a pre-determined cut-off frequency, the convolution annihilation relationship (3) can be expressed in matrix form as

$$\mathsf{Toep}(\hat{\rho})\, \boldsymbol{h} = 0. \tag{4}$$

Here $\mathsf{Toep}(\hat{\rho}) \in \mathbb{C}^{M \times N}$ denotes a rectangular Toeplitz matrix built from $\hat{\rho}[k]$, and $\boldsymbol{h} \in \mathbb{C}^N$ is a vector of the Fourier coefficients of $\mu$, zero-padded if necessary to have length $N$. This shows $\boldsymbol{h}$ is a non-trivial nullspace vector for $\mathsf{Toep}(\hat{\rho})$, i.e., $\mathsf{Toep}(\hat{\rho})$ is rank deficient. Notice also that any multiple of $\mu(x)$ by a phase factor, $\gamma(x) = \mu(x)e^{j2\pi kx}$, $k \in \mathbb{Z}$, will also satisfy the above annihilation relationship, i.e., $\mathsf{Toep}(\hat{\rho})\, \boldsymbol{h}' = 0$ where $\boldsymbol{h}'$ is a vector of the Fourier coefficients of $\gamma$. This implies that $\mathsf{Toep}(\hat{\rho})$ has a large nullspace, hence is a low-rank matrix. In fact, one can show $\text{rank}[\mathsf{Toep}(\hat{\rho})] = r$, which establishes a one-to-one correspondence between the sparsity of the signal (1) and the rank of a Toeplitz matrix built from its Fourier coefficients.

When we only have samples of $\hat{\rho}[k]$ for $k \in \Omega$ where $\Omega \subset \mathbb{Z}$ is a sampling set of arbitrary locations, we can use the low-rank property of $\mathsf{Toep}(\hat{\rho})$ to recover $\hat{\rho}[k]$ for all $|k| \leq K$ as the solution to the following matrix completion problem:

$$\min_{\hat{\phi}} \; \text{rank}[\mathsf{Toep}(\hat{\phi})] \text{ subject to } \hat{\phi}[k] = \hat{\rho}[k], \quad \forall k \in \Omega \tag{5}$$

Recovery guarantees for a convex relaxation to (5) were studied in [3] assuming the sampling locations $\Omega$ are drawn uniformly at random from the set $\{k : |k| \leq K\}$ for some fixed cut-off $K$. Multi-dimensional generalizations of the model (1) and the recovery program (5) were also considered in [3], and have been adapted to super-resolution imaging context as well [30].

Additionally, SLRA models are central to several MRI reconstruction tasks. An important example is in parallel MRI, where Fourier data of the image is collected simultaneously at multiple receive coils having different spatial sensitivity profiles. Several auto-calibration techniques in parallel MRI, such as GRAPPA [12], ESPIRiT [14], PRUNO [15], and related techniques [13], model the coil sensitivities as smooth functions in spatial domain, which translates to the low-rank property of a convolutional structured matrix built from data obtained at a uniformly sampled calibration region in Fourier domain; see, e.g., [13]. Estimates of the coil sensitivity maps are obtained from the nullspace of this matrix, and then used to reconstruct missing Fourier data via linear prediction. The extension to recovery of parallel MRI without a fully sampled calibration region was formulated as a structured low-rank matrix recovery problem in the SAKE framework [2]. Extensions of this work to the recovery of single coil data was proposed in [4], by modeling the image as having sparse support or smoothly varying phase. This was later incorporated into a multi-coil framework in [10]. Further generalizations of this approach were proposed in [6], which employed transform sparsity of the image, using wavelets and other derivative-like operators.

In the remainder of this paper, we focus on the SLRA model that motivated this work—the continuous domain piecewise smooth image model introduced in [8], [9], [24]. In [8] we showed this class of signals satisfies a Fourier domain convolution annihilation relationship provided the signal discontinuities are localized to the zero-set of a smooth bandlimited function, as in the finite-rate-of-innovation curves model [38]. This work can also be viewed as an extension of the SLRA models investigated in [39]–[42] in the context of super-resolution MRI, which was formulated for 1-D piecewise polynomial signals. Our work generalizes this model to the recovery of 2-D piecewise polynomial images with discontinuities supported on curves. For example, suppose $f(x, y)$ is a piecewise constant function in 2-D with discontinuities contained in the zero-set $\{(x, y) \in [0,1]^2 : \mu(x, y) = 0\}$ where $\mu$ is bandlimited. Analogous to (2), one can show the partial derivatives of $f$ are annihilated by multiplication with $\mu$ in spatial domain:

$$(\partial_x f)\mu = 0 \text{ and } (\partial_y f)\mu = 0 \tag{6}$$

in the distributional sense; see Figure 1 for an illustration. This translates to the Fourier domain convolution annihilation relationship:

$$(\widehat{\partial_x f} * \hat{\mu})[k_x, k_y] = 0, \tag{7}$$

$$(\widehat{\partial_y f} * \hat{\mu})[k_x, k_y] = 0 \quad \text{for all} \quad (k_x, k_y) \in \mathbb{Z}^2, \tag{8}$$

Here the expressions $\widehat{\partial_x f}$ and $\widehat{\partial_y f}$ are computed in Fourier domain by the weightings

$$\widehat{\partial_x f}[k_x, k_y] = j2\pi k_x \hat{f}[k_x, k_y],$$
$$\widehat{\partial_y f}[k_x, k_y] = j2\pi k_y \hat{f}[k_x, k_y] \quad \text{for all} \quad (k_x, k_y) \in \mathbb{Z}^2,$$

We can represent the convolution relations (7) and (8) in matrix notation as:

$$\mathcal{G}(\hat{f})\, \boldsymbol{h} = \begin{pmatrix} \mathsf{Toep}_2(\widehat{\partial_x f}) \\ \mathsf{Toep}_2(\widehat{\partial_y f}) \end{pmatrix} \boldsymbol{h} = 0 \tag{9}$$

where $\mathsf{Toep}_2(\hat{g})$ denotes the Toeplitz-like structured matrix built from any 2-D array of coefficients $\hat{g}$ representing 2-D linear convolution with $\hat{g}$. Figure 1 illustrates the construction of $\mathcal{G}(\hat{f})$, which we call the *gradient weighted matrix lifting*. Similar to the 1-D setting we can show that $\mathcal{G}(\hat{f})$ is a low-rank matrix. In [9] we proved that under certain geometric restrictions on the edge-set of $f$, there exists a unique bandlimited annihilating function $\mu_0$ whose shifts in Fourier domain span the nullspace of $\mathcal{G}(\hat{f})$. In particular, if $\mathcal{G}(\hat{f})$ is built assuming $\boldsymbol{h} \in \mathbb{C}^N$ is an array of size $\Lambda \subset \mathbb{Z}^2$, $|\Lambda| = N$, then for sufficiently large $M$ one has

$$\text{rank}(\mathcal{G}(\hat{f})) = N - S, \tag{10}$$

where $S$ is the number of integer shifts of $\widehat{\mu_0}$ inside the index set $\Lambda$ (see Prop. 6 in [9]). This establishes a correspondence between the rank of $\mathcal{G}(\hat{f})$ and the complexity of the edge-set of $f$, as measured by the bandwidth of the annihilating function $\mu_0$.

In [8] we proposed to recover the Fourier coefficients of a piecewise constant image using a structured low-rank matrix completion formulation similar to (5). This approach has several applications to under-sampled MRI reconstruction, since MRI measurements are accurately modeled as the multi-dimensional Fourier coefficients of the underlying image. See also [4]–[7] for similar Fourier domain SLRA models for MRI reconstruction based on alternative image models.

### *B. Problem Formulation*

A SLRA model supposes the data $\boldsymbol{x}_0$ to be recovered is such that a structured matrix $\mathcal{T}(\boldsymbol{x}_0)$ constructed from $\boldsymbol{x}_0$ is low-rank. In many SLRA models, including those for robust spectral estimation [3] and the super-resolution piecewise constant image model [9] introduced above, it is also known that $\mathcal{T}(\boldsymbol{x}_0)$ is the unique rank minimizer subject to certain data constraints.

This suggests we can attempt to recover $\boldsymbol{x}_0$ from its linear measurements $\boldsymbol{A}\boldsymbol{x}_0 = \boldsymbol{b}$ by solving the following rank minimization problem:

$$\min_{\boldsymbol{x}} \ \text{rank}[\mathcal{T}(\boldsymbol{x})] \ \text{ subject to } \ \boldsymbol{A}\boldsymbol{x} = \boldsymbol{b}. \tag{11}$$

However, it is well-known that (11) is NP-hard in general [43]. Many authors have investigated tractable methods to obtain exact or approximate minimizers to the structured low-rank matrix recovery problem (11), including non-convex methods based on alternating projections (also known as Cadzow's method)

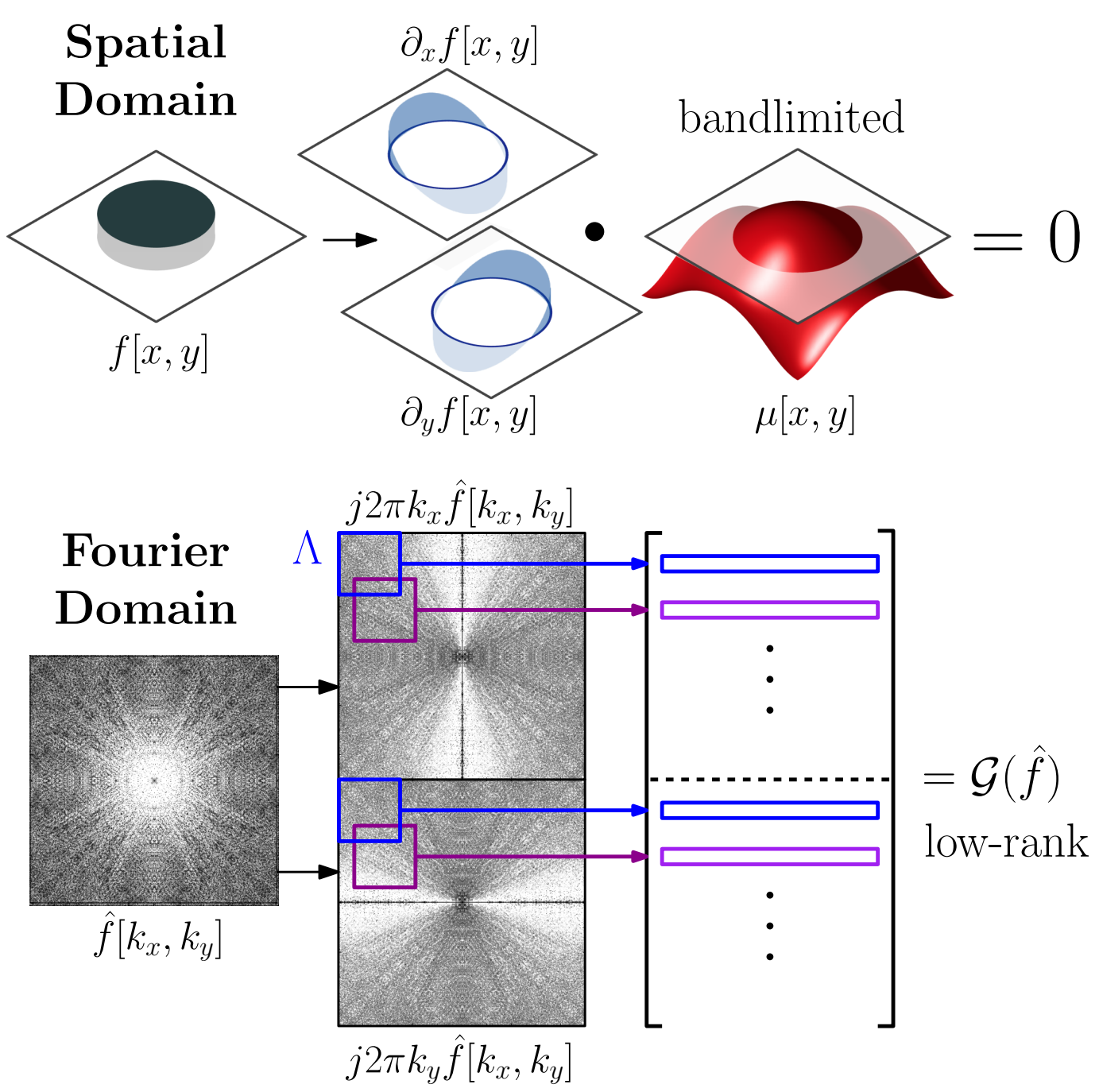

Fig. 1. 2-D annihilation relationship for piecewise constant images (top) and construction of gradient weighted lifting (bottom). If the piecewise constant image $f$ has edges supported in the zero-set of a bandlimited function $\mu$, then the partial derivatives of $f$ are annihilated by multiplication with $\mu$. This can be expressed equivalently in Fourier domain as the low-rank property of a structured matrix $\mathcal{G}(\hat{f})$ built from the Fourier coefficients $\hat{f}$. The matrix $\mathcal{G}(\hat{f})$ is constructed by weighting $\hat{f}$ to create two arrays of the simulated Fourier derivatives $\widehat{\partial_x f}$ and $\widehat{\partial_y f}$, then extracting each patch of size $\Lambda$ by a sliding window operation.

[2], [4], [44], convex relaxation methods [7], [45], or local optimization techniques [46]. We discuss several of these methods in the Supplementary Materials.

We choose to formulate the structured low-rank matrix recovery problem as a family of relaxations to (11):

$$\min_{\boldsymbol{x}} \|\mathcal{T}(\boldsymbol{x})\|_p^p \text{ subject to } \boldsymbol{A}\boldsymbol{x} = \boldsymbol{b}, \tag{12}$$

where $\|\cdot\|_p$ denotes the family of Schatten-$p$ quasi-norms with $0 < p \leq 1$, defined for an arbitrary matrix $\boldsymbol{X}$ by

$$\|\boldsymbol{X}\|_p := \left(\sum_i \sigma_i(\boldsymbol{X})^p\right)^{\frac{1}{p}}, \tag{13}$$

where $\sigma_i(\boldsymbol{X})$ are the singular values of $\boldsymbol{X}$. In the case $p = 1$, this penalty is in fact a norm, coinciding nuclear norm of a matrix. Equivalently, the Schatten-$p$ quasi-norms can be defined in terms of the following trace formula: [1]

$$\|\boldsymbol{X}\|_p = Tr[(\boldsymbol{X}^*\boldsymbol{X})^{\frac{p}{2}}]^{\frac{1}{p}}. \tag{14}$$

We also define the penalty

$$\|\boldsymbol{X}\|_0^0 := \sum_i \log(\sigma_i(\boldsymbol{X})) = \frac{1}{2}\log\det(\boldsymbol{X}^*\boldsymbol{X}) \tag{15}$$

where $\log$ denotes the natural logarithm, which can be viewed as the limiting case[2] of (13) as $p \to 0$. Note that the penalty $\|\cdot\|_1$ is convex, but the penalty $\|\cdot\|_p^p$ is non-convex for $0 \leq p < 1$.

[1]The fractional $q$-th power of a positive semi-definite matrix $\boldsymbol{Y} \in \mathbb{C}^{n\times n}$ is defined as $\boldsymbol{Y}^q = \sum_i^n \lambda_i^q \boldsymbol{v}_i\boldsymbol{v}_i^*$, where $\{\boldsymbol{v}_i\}_{i=1}^n$ is a orthonormal set of eigenvectors for $\boldsymbol{Y}$ with associated non-negative eigenvalues $\{\lambda_i\}_{i=1}^n$.

[2]More precisely, for any scalar $t > 0$ we have $\lim_{p\to 0} \frac{t^p - 1}{p} = \log(t)$. Therefore if $\boldsymbol{X} \in \mathbb{C}^{M\times N}$ has no zero singular values, $\lim_{p\to 0} \frac{1}{p}\|\boldsymbol{X}\|_p^p - \frac{N}{p} = \lim_{p\to 0}\sum_{i=1}^N \frac{\sigma_i(\boldsymbol{X})^p - 1}{p} = \sum_{i=1}^N \log(\sigma_i(\boldsymbol{X}))$

In the case where the measurements are corrupted by noise, i.e., $\boldsymbol{b} = \boldsymbol{A}\boldsymbol{x}_0 + \boldsymbol{n}$, where $\boldsymbol{n}$ is assumed to be a vector of i.i.d. complex white Gaussian noise, we relax the equality constraint in (12) by incorporating a data fidelity term into the objective:

$$\min_{\boldsymbol{x}} \|\boldsymbol{A}\boldsymbol{x} - \boldsymbol{b}\|^2 + \lambda \|\mathcal{T}(\boldsymbol{x})\|_p^p. \tag{16}$$

Here $\lambda > 0$ is a tunable regularization parameter balancing data fidelity and the degree to which $\mathcal{T}(\boldsymbol{x})$ is assumed to be low-rank.

### C. Convolutional structured matrix lifting model

By a *matrix lifting* we mean any linear operator $\mathcal{T}$ mapping a vector $\boldsymbol{x} \in \mathbb{C}^m$ to a "lifted matrix" $\mathcal{T}(\boldsymbol{x}) \in \mathbb{C}^{M\times N}$. In this work we focus on a general class of matrix liftings having a similar convolution structure to (4) and (9). For subsets $\Delta$, $\Lambda$ and $\Gamma$ of $\mathbb{Z}^2$, to be defined in the sequel, let $\mathsf{Toep}_d(\boldsymbol{y})$ be the matrix representing linear convolution with the $d$-dimensional array $\boldsymbol{y} = \{y[\boldsymbol{k}] : \boldsymbol{k} \in \Delta \subset \mathbb{Z}^d\}$, such that

$$\mathsf{Toep}_d(\boldsymbol{y})\boldsymbol{h} = \boldsymbol{y} * \boldsymbol{h}, \tag{17}$$

where $\boldsymbol{h} = \{h[\boldsymbol{k}] : \boldsymbol{k} \in \Lambda \subset \mathbb{Z}^d\}$ and $\boldsymbol{y} * \boldsymbol{h}$ is the array defined by

$$(\boldsymbol{y} * \boldsymbol{h})[\boldsymbol{k}] = \sum_{\boldsymbol{\ell}\in\Lambda} y[\boldsymbol{k} - \boldsymbol{\ell}]h[\boldsymbol{\ell}], \quad \boldsymbol{k} \in \Gamma \subset \mathbb{Z}^d. \tag{18}$$

Here the convolution is restricted to the set of "valid" indices $\Gamma \subset \mathbb{Z}^d$ satisfying $\boldsymbol{k} \in \Gamma$ only if $\boldsymbol{k} - \boldsymbol{\ell} \in \Delta$ for all $\boldsymbol{\ell} \in \Lambda$, i.e., the set of indices for which the sum in (18) is well-defined. We call $\boldsymbol{h}$ a *filter*, and $\Lambda$ the *filter support*, $\Delta$ the *data support*, and $\Gamma$ the set of *valid indices*; see Figure 2 for an illustration in dimension $d = 2$. When the index sets $\Lambda$ and $\Gamma$ are rectangular, the type of matrix structure exhibited by $\mathsf{Toep}_d(\boldsymbol{y})$ has been called *multi-level Toeplitz* [47], and we adopt this term here as well. For example, $\mathsf{Toep}_2(\boldsymbol{y})$ is a block Toeplitz matrix with Toeplitz blocks, meaning the blocks are arranged in a Toeplitz pattern and each block is itself Toeplitz.

In this work we consider structured matrix liftings that have a vertical block structure where each block is multi-level Toeplitz:

$$\mathcal{T}(\boldsymbol{x}) = \begin{pmatrix} \mathcal{T}_1(\boldsymbol{x}) \\ \vdots \\ \mathcal{T}_K(\boldsymbol{x}) \end{pmatrix} = \begin{pmatrix} \mathsf{Toep}_d(\boldsymbol{M}_1\boldsymbol{x}) \\ \vdots \\ \mathsf{Toep}_d(\boldsymbol{M}_K\boldsymbol{x}) \end{pmatrix} \tag{19}$$

Here each $\boldsymbol{M}_j$, $1 \le j \le K$ denotes some linear transformation. Typical choices of the $\boldsymbol{M}_j$ include the identity, a reshaping operator, or element-wise multiplication by a set of weights. For example, the gradient weighted lifting (9) has two blocks where $\boldsymbol{M}_1$ and $\boldsymbol{M}_2$ are diagonal matrices representing weightings by Fourier derivatives $j2\pi k_x$ and $j2\pi k_y$, respectively. Different weighting schemes have also been have been considered in related formulations [5], [8], [31].

If a lifted matrix $\mathcal{T}(\boldsymbol{x})$ in the form (19) is rank deficient, then every non-trivial vector $\boldsymbol{h} \in \text{null}\ \mathcal{T}(\boldsymbol{x})$ can be interpreted as an *annihilating filter* for each $\boldsymbol{y}_j = \boldsymbol{M}_j\boldsymbol{x}$, $1 \le j \le K$, in the sense that

$$(\boldsymbol{y}_j * \boldsymbol{h})[\boldsymbol{k}] = 0 \ \text{ for all } \ \boldsymbol{k} \in \Gamma. \tag{20}$$

In other words, if the structured matrix (19) is low-rank, it means there exists a large collection of linearly independent annihilating filters for the data defining the matrix. This generalizes the annihilating filter formulation that is central to finite-rate-of-innovation modeling [48].

Finally, we note that any (multi-level) Hankel matrix can be rewritten as a (multi-level) Toeplitz matrix through a permutation of its rows and columns, which has no effect on the rank of the matrix. In this way we can also incorporate (block) multi-level Hankel matrix liftings into the model (19), such as those proposed in [3], [6].

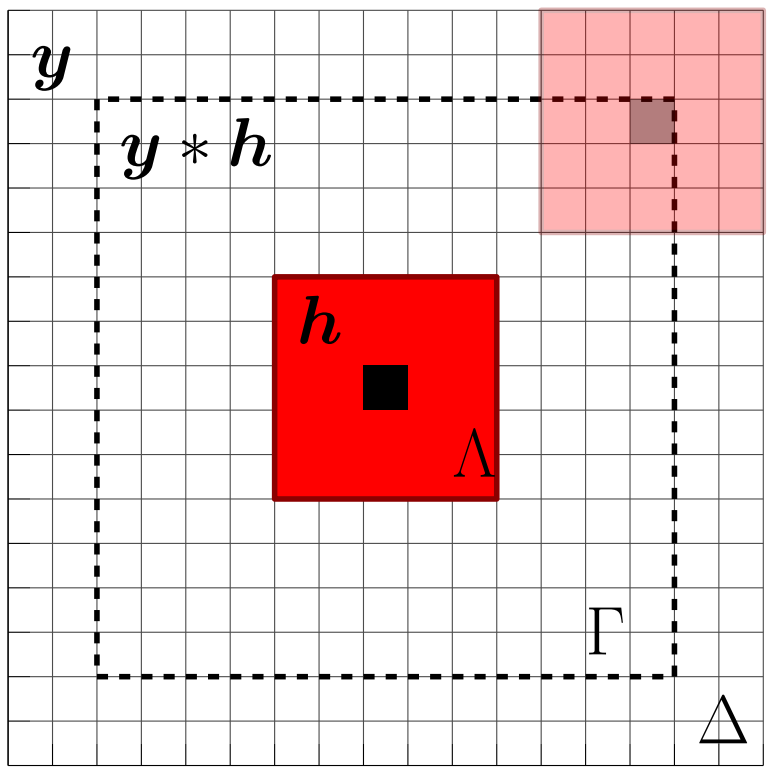

Fig. 2. Illustration of index sets used in construction of multi-level Toeplitz matrix liftings (17) in two-dimensions ($d = 2$). Here $\Delta$ is the support of the data $\boldsymbol{y}$, $\Lambda$ (in red) is the support of the filter $\boldsymbol{h}$ with index $(0, 0)$ in black, and $\Gamma$ (interior of dashed line) represents the valid set where the linear convolution $\boldsymbol{y} * \boldsymbol{h}$ is well-defined.

## III. Iteratively reweighted least-squares algorithms for structured low-rank matrix recovery

We adapt an iteratively reweighted least-squares (IRLS) approach to minimizing (12) or (16), originally proposed in [20], [21] in the low-rank matrix completion setting; see also [49] for an alternative matrix factorization approach to Schatten-$p$ norm minimization. The IRLS approach is motivated by the observation that the Schatten-$p$ quasi-norm of any matrix $\boldsymbol{X}$ may be re-expressed as a weighted Frobenius norm, where the weights depend on $\boldsymbol{X}$:

$$\|\boldsymbol{X}\|_p^p = Tr[(\boldsymbol{X}^*\boldsymbol{X})\underbrace{(\boldsymbol{X}^*\boldsymbol{X})^{\frac{p}{2}-1}}_{\boldsymbol{H}}] = \|\boldsymbol{X}\sqrt{\boldsymbol{H}}\|_F^2, \tag{21}$$

for all $0 < p \leq 1$, provided $\boldsymbol{X}$ has no zero singular values so that $\boldsymbol{H}$ is well-defined. Here $\sqrt{\boldsymbol{H}}$ denotes any matrix[3] satisfying $\boldsymbol{H} = (\sqrt{\boldsymbol{H}})(\sqrt{\boldsymbol{H}})^*$. This suggests an iterative algorithm for minimizing the Schatten-$p$ quasi-norm that alternates between updating a weight matrix $\boldsymbol{H}$ and solving a weighted least-squares problem with $\boldsymbol{H}$ fixed. In particular, substituting $\boldsymbol{X} = \mathcal{T}(\boldsymbol{x})$ in the IRLS-$p$ algorithm presented in [21] gives the following iterative scheme for solving (16):

$$\boldsymbol{H}_n = [\mathcal{T}(\boldsymbol{x}^{(n-1)})^*\mathcal{T}(\boldsymbol{x}^{(n-1)}) + \epsilon_n]^{\frac{p}{2}-1} \tag{22}$$

$$\boldsymbol{x}^{(n)} = \arg\min_{\boldsymbol{x}} \|\boldsymbol{A}\boldsymbol{x} - \boldsymbol{b}\|_2^2 + \lambda C_p \|\mathcal{T}(\boldsymbol{x})\sqrt{\boldsymbol{H}_n}\|_F^2. \tag{23}$$

where we set $C_p = \frac{p}{2}$ if $0 < p \leq 1$ and $C_0 = \frac{1}{2}$, and $\epsilon_n > 0$ is an iteration dependent smoothing parameter that is typically decreased in each iteration. In order to ensure long-run stability of the algorithm, we follow the approach in [21], and decrease $\epsilon_n$ until it reaches some pre-determined minimum value $\epsilon_{\text{min}}$, after which it is remained fixed at this value, rather than shrinking $\epsilon_n$ to zero. Pseudo-code summarizing this approach is given in Algorithm 1. In the Supplementary Materials we show that for a fixed $\epsilon$ algorithm 1 can be derived as a majorization-minimization (MM) algorithm [50] for a cost function similar to (16) but with a *smoothed* version of the Schatten-$p$; the MM interpretation informs our choice of the constant $C_p$ in (23). Due to properties of MM algorithms, this ensures that after $\epsilon_n = \epsilon_{\text{min}}$, the iterates monotonically decrease the modified cost function. Moreover, in the convex case $p = 1$ the iterates are guaranteed to converge to the global minimum of the cost function. In the non-convex case $p \in [0, 1)$, convergence to the global optimality cannot be ensured, but the iterates are still guaranteed to converge to a stationary point

[3]One choice of $\sqrt{\boldsymbol{H}}$ is the standard matrix square root $\boldsymbol{H}^{\frac{1}{2}} = \sum_i \lambda_i^{\frac{1}{2}} \boldsymbol{v}_i \boldsymbol{v}_i^*$ where $(\lambda_i, \boldsymbol{v}_i)$ denotes an eigen-pair, but we will show later it is computationally advantageous to use a different choice of square root in our algorithm.

[51]. We refer the reader to [21], [22] for a more detailed convergence analysis of the IRLS-$p$ algorithm in the context of low-rank matrix completion, noting that many of these results can be translated to the structured matrix setting as well.

---

**Algorithm 1:** IRLS-$p$ algorithm for structured low-rank matrix recovery

---

Initialize $\boldsymbol{x}^{(0)}$ and choose $\epsilon_0, \epsilon_{\min} > 0$;
**for** $n = 1$ *to* $I_{max}$ **do**
  $\boldsymbol{H}_n = [\mathcal{T}(\boldsymbol{x}^{(n-1)})^*\mathcal{T}(\boldsymbol{x}^{(n-1)}) + \epsilon_{n-1}\boldsymbol{I}]^{\frac{p}{2}-1}$;
  $\boldsymbol{x}^{(n)} = \arg\min_{\boldsymbol{x}} \|\boldsymbol{A}\boldsymbol{x} - \boldsymbol{b}\|_2^2 + \lambda C_p\|\mathcal{T}(\boldsymbol{x})\sqrt{\boldsymbol{H}_n}\|_F^2$;
  Choose $\epsilon_n$ such that $0 < \epsilon_{\min} \leq \epsilon_n \leq \epsilon_{n-1}$;
**end**

---

### A. Challenges with the direct IRLS approach

While the IRLS-$p$ algorithm has shown several advantages in the low-rank matrix completion setting [20]–[22], its direct adaptation to the structured matrix liftings considered in this work is computationally prohibitive for large-scale problems. To see why, first note that the weight matrix update requires computing an inverse power of the Gram matrix $\mathcal{T}(\boldsymbol{x})^*\mathcal{T}(\boldsymbol{x}) \in \mathbb{C}^{N\times N}$. Computing this Gram matrix directly by matrix multiplication will be costly when the inner dimension $M \gg N$, and inverting the resulting matrix can be unstable due to round-off errors. Instead, a more stable and efficient method to compute the inverse is by an SVD of the $M \times N$ matrix $\mathcal{T}(\boldsymbol{x})$, which requires $O(MN^2)$ flops to compute. In the case that $\mathcal{T}(\boldsymbol{x})$ is approximately low-rank, one can instead compute only the top $r$ singular values and singular vectors either using deterministic methods, such as Lanczos bidiagonalization [52], or by randomized methods [53], at a reduced cost of $O(MNr)$ flops. However, these approaches are still dominated by a computational cost that is linear in $M$, and will be prohibitively slow when $M$ is large or when the matrix $\mathcal{T}(\boldsymbol{x})$ is not sufficiently low-rank.

The IRLS-$p$ algorithm additionally requires solving the weighted least-squares problem (23) at each iteration. To give an idea of the costs involved in solving (23) via an iterative solver, consider the case where the matrix lifting has the form $\mathcal{T}(\boldsymbol{x}) = \mathsf{Toep}_d(\boldsymbol{x})$, i.e., $\mathcal{T}(\boldsymbol{x})\boldsymbol{h} = \boldsymbol{x} * \boldsymbol{h}$. In this case, standard iterative methods for solving (23), such as the CG or LSQR algorithm [54], will require $O(N)$ multi-dimensional FFTs per iteration, where $N = |\Lambda|$ is the total number of filter coefficients. However, $N$ can be on the order of hundreds or thousands for the problems considered in this work. Therefore, standard methods for solving the least-squares problem (23) will also be prohibitively costly for these large-scale problems.

## IV. Proposed GIRAF algorithm

As observed in the previous section, the direct IRLS algorithm will not scale well to large problem instances because of two challenges: 1) computing the weight matrix update requires a large-scale SVD, and 2) computing a solution to the least-squares problem requires a prohibitive number of FFTs. In this section we propose a novel approximation of the problem formulation (16) to overcome these difficulties. The main idea is to *approximate the structured matrix lifting* $\mathcal{T}$ in a systematic way such that the complexity of the resulting IRLS subproblems simplify, while preserving the rank structure of the lifting as best as possible.

### A. Half-circulant approximation of a Toeplitz matrix

Our approximation is based on the observation that every multi-level Toeplitz matrix can be embedded in a larger multi-level circulant matrix. Specifically, the $d$-level Toeplitz matrix $\mathsf{Toep}_d(\boldsymbol{y})$ built with filter support $\Lambda \subset \mathbb{Z}^d$ and valid index set $\Gamma \subset \mathbb{Z}^d$ can always be expressed as:

$$\mathsf{Toep}_d(\boldsymbol{y}) = \boldsymbol{P}_\Gamma\, \boldsymbol{C}(\boldsymbol{y})\, \boldsymbol{P}_\Lambda^* \tag{24}$$

Here $\boldsymbol{C}(\boldsymbol{y}) \in \mathbb{C}^{L\times L}$ is a matrix representing convolution with the array $\boldsymbol{y}$, $\boldsymbol{P}_\Gamma \in \mathbb{C}^{M\times L}$ is a matrix representing restriction to valid set $\Gamma$, and $\boldsymbol{P}_\Lambda^* \in \mathbb{C}^{L\times N}$ represents a zero-padding outside the filter support $\Lambda$. The matrices $\boldsymbol{P}_\Gamma$ and $\boldsymbol{P}_\Lambda^*$ act as row-restriction and column-restriction operators, respectively; see Figure (3) for an illustration. Because of the restriction matrix $\boldsymbol{P}_\Gamma$, we can assume $\boldsymbol{C}(\boldsymbol{y})$ represents a multi-dimensional *circular convolution*, i.e., $\boldsymbol{C}(\boldsymbol{y})$ is multi-level circulant, provided the convolution takes place on a sufficiently large rectangular grid to avoid wrap-around boundary effects. In particular, the rectangular circular convolution grid should be at least as large as the data support $\Delta \subset \mathbb{Z}^d$. Without loss of generality, from now on we assume the data support $\Delta$ is rectangular so that it coincides with the circular convolution grid.

To simplify subsequent derivations, in this section we assume the structured matrix lifting $\mathcal{T}(\boldsymbol{x})$ consists of a single multi-level Toeplitz block: $\mathcal{T}(\boldsymbol{x}) = \mathsf{Toep}_d(\boldsymbol{M}\boldsymbol{x})$. According to (24), $\mathcal{T}(\boldsymbol{x})$ can be written as

$$\mathcal{T}(\boldsymbol{x}) = \boldsymbol{P}_\Gamma\, \boldsymbol{C}(\boldsymbol{M}\boldsymbol{x})\, \boldsymbol{P}_\Lambda^*. \tag{25}$$

We propose approximating the structured matrix lifting $\mathcal{T}(\boldsymbol{x})$ with the surrogate lifting $\breve{\mathcal{T}}(\boldsymbol{x})$ defined by

$$\breve{\mathcal{T}}(\boldsymbol{x}) = \boldsymbol{C}(\boldsymbol{M}\boldsymbol{x})\boldsymbol{P}_\Lambda^* \tag{26}$$

i.e., we omit the left-most row restriction operator $\boldsymbol{P}_\Gamma$ from $\mathcal{T}(\boldsymbol{x})$ so that $\breve{\mathcal{T}}(\boldsymbol{x})$ is the full vertical section of the multi-level circulant matrix $\boldsymbol{C}(\boldsymbol{M}\boldsymbol{x})$. In general, we say a matrix $\boldsymbol{X} \in \mathbb{C}^{M\times N}$, $N \leq M$, is (multi-level) *half-circulant* if it can be written as $\boldsymbol{X} = \boldsymbol{C}\boldsymbol{P}^*$, where $\boldsymbol{C} \in \mathbb{C}^{M\times M}$ is (multi-level) circulant and $\boldsymbol{P} \in \mathbb{C}^{N\times M}$ is a restriction matrix, i.e., $\boldsymbol{X}$ is obtained by selecting $N$ full columns from $\boldsymbol{C}$.

Inserting the approximation (26) into (16), we propose solving

$$\min_x \|\boldsymbol{A}\boldsymbol{x} - \boldsymbol{b}\|^2 + \lambda \|\breve{\mathcal{T}}(\boldsymbol{x})\|_p^p. \tag{27}$$

using the IRLS-$p$ algorithm. In other words, rather than penalizing the Schatten-$p$ norm of the exact multi-level Toeplitz lifted matrix, we penaltize its multi-level half-circulant approximation instead. We can justify this approach by considering the effect the half-circulant approximation (26) has on the singular values of the lifting. Recall that $\boldsymbol{C}(\boldsymbol{y})$ denotes circular convolution with the array $\boldsymbol{y}$ on a rectangular grid $\Delta \subset \mathbb{Z}^d$ containing the index set $\Gamma$. Hence, after a rearrangement of rows, we can always write $\mathcal{T}(\boldsymbol{x})$ submatrix of $\breve{\mathcal{T}}(\boldsymbol{x})$:

$$\breve{\mathcal{T}}(\boldsymbol{x}) = \begin{bmatrix} \mathcal{T}(\boldsymbol{x}) \\ \mathcal{B}(\boldsymbol{x}) \end{bmatrix} = \begin{bmatrix} \boldsymbol{P}_\Gamma \boldsymbol{C}(\boldsymbol{M}\boldsymbol{x})\boldsymbol{P}_\Lambda^* \\ \boldsymbol{P}_{\Gamma^C} \boldsymbol{C}(\boldsymbol{M}\boldsymbol{x})\boldsymbol{P}_\Lambda^* \end{bmatrix}, \tag{28}$$

where $\Gamma^C$ is the set complement of $\Gamma$ inside the circular convolution grid $\Delta$. In other words, $\breve{\mathcal{T}}(\boldsymbol{x})$ is the unique matrix obtained by augmenting the rows of $\mathcal{T}(\boldsymbol{x})$ to make it (multi-level) half-circulant.

As a consequence of the embedding (28), the singular values of $\mathcal{T}(\boldsymbol{x})$ are bounded by those of $\breve{\mathcal{T}}(\boldsymbol{x})$:

$$\sigma_i(\mathcal{T}(\boldsymbol{x})) \leq \sigma_i(\breve{\mathcal{T}}(\boldsymbol{x})) \quad \text{for} \quad i = 1, \ldots, N. \tag{29}$$

(see Corollary 3.1.3 in [55]). Hence, the Schatten-$p$ quasi-norm of $\mathcal{T}(\boldsymbol{x})$ is also bounded by that of $\breve{\mathcal{T}}(\boldsymbol{x})$:

$$\|\mathcal{T}(\boldsymbol{x})\|_p^p \leq \|\breve{\mathcal{T}}(\boldsymbol{x})\|_p^p \;\; \text{for all} \;\; 0 \leq p \leq 1. \tag{30}$$

Moreover, if in minimizing $\|\breve{\mathcal{T}}(\boldsymbol{x})\|_p^p$ we obtain a low-rank matrix $\breve{\mathcal{T}}(\boldsymbol{x}')$ with rank $r$, the singular value bounds (29) show that $\mathcal{T}(\boldsymbol{x}')$ is also low-rank with rank $\leq r$. This shows that it is reasonable to use $\|\breve{\mathcal{T}}(\boldsymbol{x})\|_p^p$ as a surrogate penalty for $\|\mathcal{T}(\boldsymbol{x})\|_p^p$.

Empirically, we find that when the data support is sufficiently large relative to the filter size, i.e., when $M \gg N$ where $\mathcal{T}(\boldsymbol{x}) \in \mathbb{C}^{M\times N}$, then using the surrogate lifting $\breve{\mathcal{T}}(\boldsymbol{x})$ results in negligible approximation errors (see Figure 5). When the data support $\Delta$ is small, we recommend solving for the data array $\boldsymbol{x}$ on a slightly larger "oversampled" grid $\Delta'$ to keep the approximation error low. Typically we take $\Delta' = \Delta + 2\Lambda$, i.e. we pad the data support with an extra margin having the size of the filter support. After solving the problem on the oversampled grid $\Delta'$, we finally restrict the solution to the desired data support $\Delta$. We study the effect of using an oversampled grid in Section V (see Figure 6).

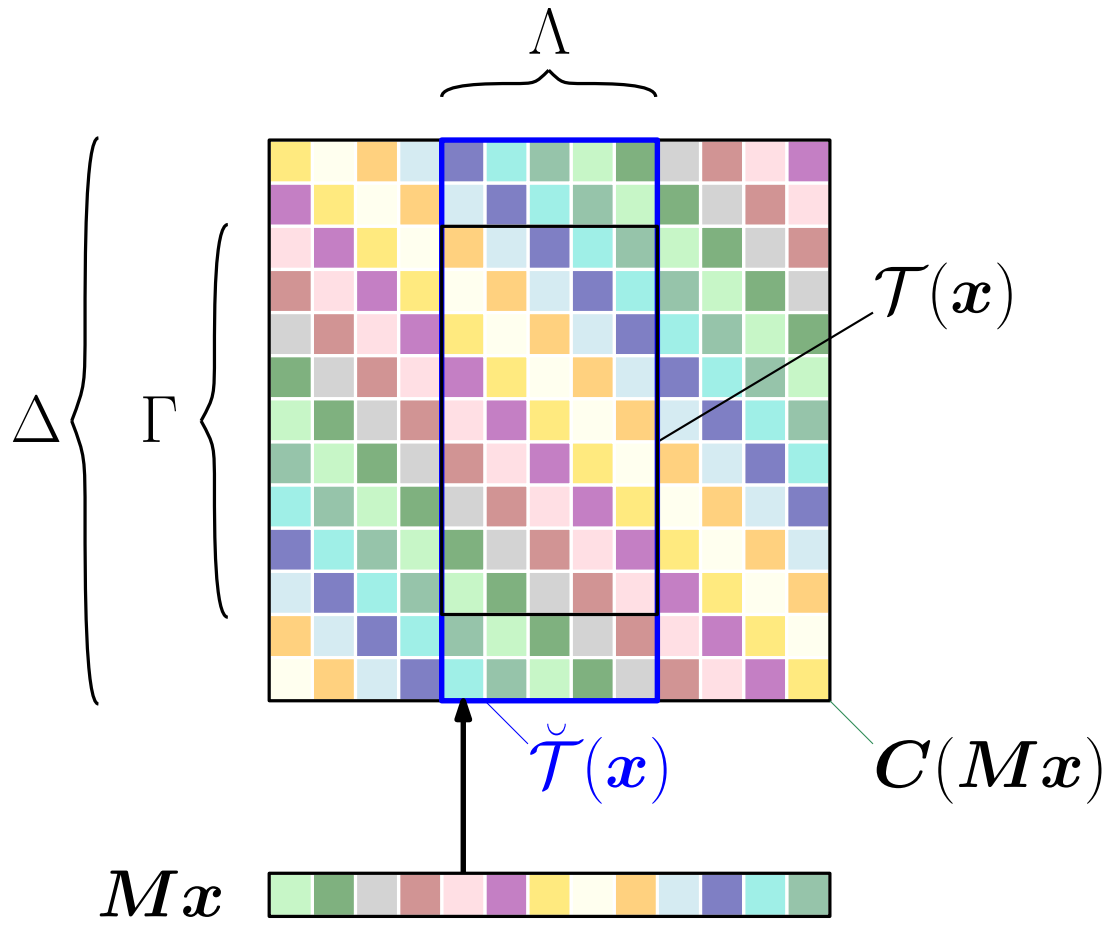

Fig. 3. *Example of Half-circulant approximation in 1-D.* We approximate the Toeplitz matrix lifting $\mathcal{T}(\boldsymbol{x})$ with the half-circulant matrix $\breve{\mathcal{T}}(\boldsymbol{x})$ obtained by adding rows to $\mathcal{T}(\boldsymbol{x})$ to make it the full vertical section of a circulant matrix $\boldsymbol{C}(\boldsymbol{M}\boldsymbol{x})$

We now show the IRLS algorithm applied to (27) results in subproblems with significantly reduced complexity due to the half-circulant structure of the approximated lifting. In particular, we will show the algorithm can be interpreted as alternating between: (1) solving for an annihilating filter for the data, and (2) solving for the data best annihilated by the filter in a least-squares sense. Due to this interpretation, we call this approach the Generic Iterative Reweighted Annihilating Filter (GIRAF) algorithm.

## *B. GIRAF least-squares subproblem*

*1) Reformulation as least-squares annihilation:* For a fixed weight matrix $\boldsymbol{H}$, the GIRAF least-squares subproblem has the form

$$\min_{\boldsymbol{x}} \|\boldsymbol{A}\boldsymbol{x}-\boldsymbol{b}\|_2^2 + \lambda C_p \|\breve{\mathcal{T}}(\boldsymbol{x})\sqrt{\boldsymbol{H}}\|_F^2. \tag{31}$$

Under the half-circulant assumption (26), we can re-express the Frobenius norm above as

$$\|\breve{\mathcal{T}}(\boldsymbol{x})\sqrt{\boldsymbol{H}}\|_F^2 = \sum_{i=1}^N \|\breve{\mathcal{T}}(\boldsymbol{x})\boldsymbol{h}_i\|_2^2 = \sum_{i=1}^N \|\boldsymbol{C}(\boldsymbol{M}\boldsymbol{x})\boldsymbol{P}_\Lambda^*\boldsymbol{h}_i\|_2^2 \tag{32}$$

where $\boldsymbol{h}_i$ is the $i$th column of $\sqrt{\boldsymbol{H}}$. Because circular convolution is commutative, we can write $\boldsymbol{C}(\boldsymbol{M}\boldsymbol{x})\boldsymbol{P}_\Lambda^*\boldsymbol{h}_i = \boldsymbol{C}_i\boldsymbol{M}\boldsymbol{x}$, where $\boldsymbol{C}_i := \boldsymbol{C}(\boldsymbol{P}_\Lambda^*\boldsymbol{h}_i)$ is a multi-level circulant matrix representing circular convolution with the zero-padded filter $\boldsymbol{h}_i$ on the rectangular grid $\Delta \subset \mathbb{Z}^d$. Therefore, $\boldsymbol{C}_i = \boldsymbol{F}\boldsymbol{D}_i\boldsymbol{F}^*$ where $\boldsymbol{F}$ is the DFT on $\Delta$, and $\boldsymbol{D}_i$ the diagonal matrix having entries given by the array $\boldsymbol{d}_i = \boldsymbol{F}^*\boldsymbol{P}_\Lambda^*\boldsymbol{h}_i$, the inverse DFT of the zero-padded filter $\boldsymbol{h}_i$. This allows us to simplify the right-hand side of (32) as

$$\begin{aligned}\sum_{i=1}^N \|\boldsymbol{C}_i\boldsymbol{M}\boldsymbol{x}\|_2^2 &= \boldsymbol{x}^*\boldsymbol{M}^*\left(\sum_{i=1}^N \boldsymbol{C}_i^*\boldsymbol{C}_i\right)\boldsymbol{M}\boldsymbol{x}\\ &= \boldsymbol{x}^*\boldsymbol{M}^*\boldsymbol{F}\left(\sum_{i=1}^N \boldsymbol{D}_i^*\boldsymbol{D}_i\right)\boldsymbol{F}^*\boldsymbol{M}\boldsymbol{x}\\ &= \|\boldsymbol{D}^{\frac{1}{2}}\boldsymbol{F}^*\boldsymbol{M}\boldsymbol{x}\|_2^2,\end{aligned}$$

where we have set $\boldsymbol{D} = \sum_{i=1}^N \boldsymbol{D}_i^*\boldsymbol{D}_i$, which is again diagonal with entries $\boldsymbol{d} \in \mathbb{C}^{|\Delta|}$ given by

$$\boldsymbol{d} = \sum_{i=1}^N |\boldsymbol{d}_i|^2, \quad \boldsymbol{d}_i = \boldsymbol{F}^*\boldsymbol{P}_\Lambda^*\boldsymbol{h}_i \quad \text{for all} \quad i = 1, ..., N, \tag{33}$$

where $|\cdot|$ is applied element-wise. Therefore, (31) transforms into the weighted least-squares problem:

$$\min_{\boldsymbol{x}} \|\boldsymbol{A}\boldsymbol{x}-\boldsymbol{b}\|_2^2 + \lambda C_p \|\boldsymbol{D}^{\frac{1}{2}}\boldsymbol{F}^*\boldsymbol{M}\boldsymbol{x}\|_2^2. \tag{34}$$

Observe that by making use of the half-circulant assumption we have effectively reduced the working dimension of the least-squares problem back down to the dimension of the original decision variable $\boldsymbol{x}$.

Computing the weights $\boldsymbol{d}$ according to the formula (33) can be costly for large-scale problems since it requires $N$ large-scale FFTs. A more efficient approach is to pre-compute the filter $\boldsymbol{h}$ defined by

$$\boldsymbol{h} = \sum_{i=1}^{N} (\widetilde{\boldsymbol{h}}_i * \boldsymbol{h}_i) \tag{35}$$

where $\widetilde{\boldsymbol{h}}_i$ is the reversed, conjugated filter defined by $\widetilde{\boldsymbol{h}}_i[\boldsymbol{k}] = \overline{\boldsymbol{h}_i[-\boldsymbol{k}]}$ for all $\boldsymbol{k} \in \Lambda$. Note that $\boldsymbol{h}$ has coefficients supported within $2\Lambda := \{\boldsymbol{k}+\boldsymbol{\ell} : \boldsymbol{k},\boldsymbol{\ell} \in \Lambda\}$, since each $\boldsymbol{h}_i$ is supported within $\Lambda$. Therefore, applying the DFT convolution theorem to (33), $\boldsymbol{d}$ can be obtained by $\boldsymbol{d} = \boldsymbol{F}^*\boldsymbol{P}_{2\Lambda}^*\boldsymbol{h}$, which after computing $\boldsymbol{h}$, requires only one large-scale FFT. We call filter $\boldsymbol{h}$ defined in (35) the *re-weighted annihilating filter* and $\boldsymbol{d}$ the *annihilation weights*.

*2) ADMM solution of least-squares annihilation:* The linear least-squares problem (34) can be readily solved with an off-the-shelf solver, such as the CG or LSQR algorithm [54]. However, the problem is poorly conditioned when the annihilation weights $\boldsymbol{d}$ are close to zero, which is expected to happen as the iterations progress. Therefore, solving (34) efficiently during the course of the GIRAF algorithm will require a robust preconditioning strategy. However, because the transformation $\boldsymbol{M}$ is not always invertible, designing an all-purpose preconditioner for (34) is challenging. Instead, we adopt an approach which allows us to solve a series of subproblems with predictably good conditioning. The approach is based on the following variable splitting:

$$\min_{\boldsymbol{x}} \|\boldsymbol{A}\boldsymbol{x}-\boldsymbol{b}\|_2^2 + \lambda C_p \|\boldsymbol{D}^{\frac{1}{2}}\boldsymbol{y}\|_2^2 \;\text{ subject to }\; \boldsymbol{F}\boldsymbol{y} = \boldsymbol{M}\boldsymbol{x}. \tag{36}$$

The equality constrained problem (36) can be efficiently solved with the alternating directions method of mulipliers (ADMM) algorithm [56], which results in the following iterative scheme:

$$\boldsymbol{y}^{(n)} = \arg\min_{\boldsymbol{y}} \left\{ \|\boldsymbol{D}^{\frac{1}{2}}\boldsymbol{y}\|_2^2 + \gamma\|\boldsymbol{y} - \boldsymbol{F}^*\boldsymbol{M}\boldsymbol{x}^{(n-1)} + \boldsymbol{q}^{(n-1)}\|_2^2 \right\}, \tag{37}$$

$$\boldsymbol{x}^{(n)} = \arg\min_{\boldsymbol{x}} \left\{ \|\boldsymbol{A}\boldsymbol{x}-\boldsymbol{b}\|_2^2 + \gamma\lambda C_p\|\boldsymbol{y}^{(n)} - \boldsymbol{F}^*\boldsymbol{M}\boldsymbol{x} + \boldsymbol{q}^{(n-1)}\|_2^2 \right\}, \tag{38}$$

$$\boldsymbol{q}^{(n)} = \boldsymbol{q}^{(n-1)} + \boldsymbol{F}\boldsymbol{y}^{(n)} - \boldsymbol{M}\boldsymbol{x}^{(n)}, \tag{39}$$

where $\boldsymbol{q}$ represents a vector of Lagrange multipliers, and $\gamma > 0$ is a fixed parameter that can be tuned to improve the conditioning of the subproblems. Subproblems (37) and (38) are both quadratic and can be solved efficiently: (37) has the exact solution

$$\boldsymbol{y}^{(n)} = (\boldsymbol{D} + \gamma\boldsymbol{I})^{-1}(\gamma\boldsymbol{F}^*[\boldsymbol{M}\boldsymbol{x}^{(n-1)} - \boldsymbol{q}]). \tag{40}$$

Since $\boldsymbol{D}+\gamma\boldsymbol{I}$ is diagonal, its inverse acts as an element-wise division. Likewise, the solution of (38) can be obtained as

$$\boldsymbol{x}^{(n)} = (\boldsymbol{A}^*\boldsymbol{A} + \gamma\lambda C_p\boldsymbol{M}^*\boldsymbol{M})^{-1}\left[\boldsymbol{A}^*\boldsymbol{b} + \gamma\lambda C_p\boldsymbol{M}^*\boldsymbol{F}(\boldsymbol{y}^{(n)} + \boldsymbol{q}^{(n-1)})\right]. \tag{41}$$

In many problems of interest, both $\boldsymbol{A}^*\boldsymbol{A}$ and $\boldsymbol{M}^*\boldsymbol{M}$ are diagonal, which means $\boldsymbol{x}^{(n)}$ is also obtained by an efficient element-wise division. In this case, aside from the element-wise operations, the computational cost of one pass of ADMM iterations (37)–(39) is three FFTs. When either $\boldsymbol{A}^*\boldsymbol{A}$ or $\boldsymbol{M}^*\boldsymbol{M}$ are not diagonal, an approximate solution to (38) can be found by a few passes of an iterative solver instead. Even with inexact updates of the $\boldsymbol{x}$-subproblem, convergence of the overall ADMM algorithm is still guaranteed under fairly broad conditions; see, e.g., [56], [57].

Update (40) suggests choosing $\gamma$ adaptively before solving the least-squares subproblem according to the magnitude of the annihilation weights $\boldsymbol{d} = \{d_j\}_{j=1}^L$ of $\boldsymbol{D}$ to ensure good conditioning and fast convergence of the ADMM scheme. We recommend choosing $\gamma = (\max_j d_j)/\delta$ for $\delta \geq 1$; we study the effect of $\delta$ in Section V.

### *C. GIRAF re-weighted annihilating filter subproblem*

*1) Construction of Gram matrix:* At each iteration, the GIRAF algorithm requires updating a weight matrix $\boldsymbol{H}$ according to:

$$\boldsymbol{H} = (\breve{\mathcal{T}}(\boldsymbol{x})^*\breve{\mathcal{T}}(\boldsymbol{x}) + \epsilon\,\boldsymbol{I})^{\frac{p}{2}-1} \tag{42}$$

Rather than computing $\boldsymbol{H}$ via a full or partial SVD of the tall $M \times N$ matrix $\breve{\mathcal{T}}(\boldsymbol{x})$, as is recommended in [21], we propose computing $\boldsymbol{H}$ via the eigendecomposition of the smaller $N \times N$ Gram matrix $\boldsymbol{G} = \breve{\mathcal{T}}(\boldsymbol{x})^*\breve{\mathcal{T}}(\boldsymbol{x})$. This is because the half-circulant approximation (26) allows us to compute $\boldsymbol{G}$ efficiently in a matrix-free manner, as we now show. From (26), we have

$$\boldsymbol{G} = \boldsymbol{P}_\Lambda \boldsymbol{C}(\boldsymbol{M}\boldsymbol{x})^*\boldsymbol{C}(\boldsymbol{M}\boldsymbol{x})\boldsymbol{P}_\Lambda^*. \tag{43}$$

The restriction matrices $\boldsymbol{P}_\Lambda$ and $\boldsymbol{P}_\Lambda^*$ in (43) extract an $N \times N$ block of $\boldsymbol{C}(\boldsymbol{M}x)^*\boldsymbol{C}(\boldsymbol{M}x)$ corresponding to the intersection of the rows and columns indexed by the filter support set $\Lambda$. Note that the product $\boldsymbol{C}(\boldsymbol{M}\boldsymbol{x})^*\boldsymbol{C}(\boldsymbol{M}\boldsymbol{x}) = \boldsymbol{C}(\boldsymbol{g})$ is again a multi-level circulant matrix whose entries are generated by the array $\boldsymbol{g} = \boldsymbol{F}|\boldsymbol{F}^*\boldsymbol{M}\boldsymbol{x}|^2$ with the operation $|\cdot|^2$ is understood element-wise. Therefore, we can build $\boldsymbol{G}$ by performing a sliding-window operation that extracts every patch of size $\Lambda$ from the array $\boldsymbol{g}$. In particular, because of the restriction matrices $\boldsymbol{P}_\Lambda$ and $\boldsymbol{P}_\Lambda^*$ in (43), we only need to consider patches coming from $\boldsymbol{g}$ restricted to the index set $2\Lambda = \{\boldsymbol{k} + \boldsymbol{\ell} : \boldsymbol{k}, \boldsymbol{\ell} \in \Lambda\}$. Aside from this sliding-window operation, the main cost in computing $\boldsymbol{G}$ is two FFTs.

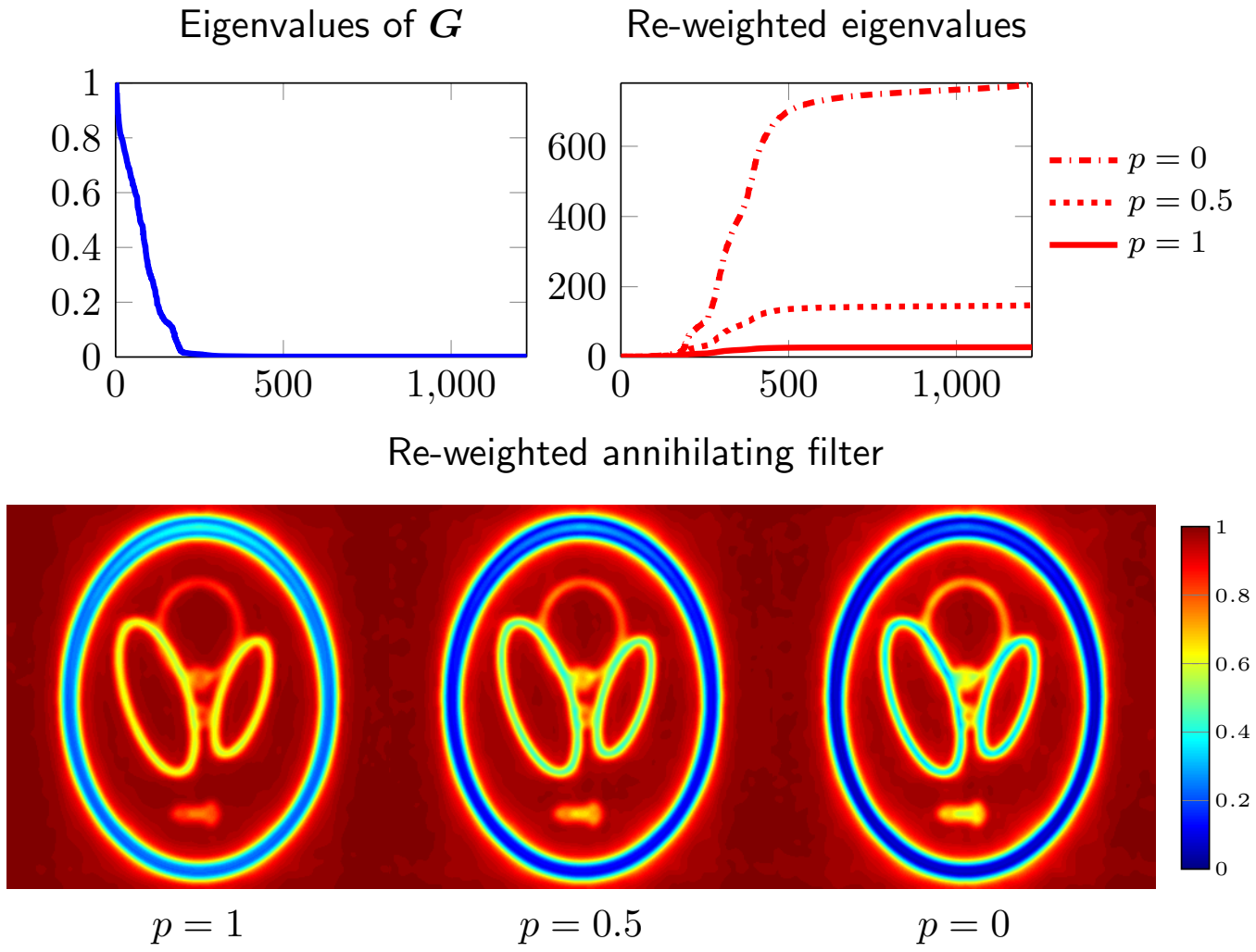

Fig. 4. Effect of Schatten-$p$ penalty on the GIRAF re-weighted annihilating filter. We show the eigenvalues of $\boldsymbol{G} = \breve{\mathcal{T}}(\boldsymbol{x})^*\breve{\mathcal{T}}(\boldsymbol{x})$ during one iteration of the GIRAF algorithm applied to the gradient weighted lifting (9) (top left), and the different re-weightings of the eigenvalues specified by different Schatten-$p$ penalties (top right). The set of annihilation weights determined by the resulting re-weighted annihilating filter given by (46) are shown in the bottom row, normalized to $[0, 1]$. Note that for smaller values of $p$ the annihilation weights are closer to zero near the edges because the re-weighted annihilating filter gives higher weight to filters in the nullspace of $\breve{\mathcal{T}}(\boldsymbol{x})$.

*2) Update of re-weighted annihilating filter:* According to section IV-B, the GIRAF least-squares subproblem can be interpreted as an annihilation of the data subject to a filter $\boldsymbol{h}$ determined by the columns $\boldsymbol{h}_1, ..., \boldsymbol{h}_N$ of $\sqrt{\boldsymbol{H}}$. We now show how to update $\boldsymbol{h}$ directly from an eigendecomposition of the Gram matrix $\boldsymbol{G} = \breve{\mathcal{T}}(\boldsymbol{x})^*\breve{\mathcal{T}}(\boldsymbol{x})$, rather than forming $\sqrt{\boldsymbol{H}}$ explicitly. Let $\boldsymbol{G} = \boldsymbol{V}\boldsymbol{\Lambda}\boldsymbol{V}^*$ where the columns $\{\boldsymbol{v}_i\}_{i=1}^N$ of $\boldsymbol{V}$ are an orthonormal basis of eigenvectors and $\boldsymbol{\Lambda}$ is diagonal matrix of the associated eigenvalues $\{\lambda_i\}_{i=1}^N$.

Then the weight matrix update (42) reduces to

$$\boldsymbol{H} = [\boldsymbol{V}(\boldsymbol{\Lambda}+\epsilon\boldsymbol{I})\boldsymbol{V}^*]^{\frac{p}{2}-1} = \boldsymbol{V}(\boldsymbol{\Lambda}+\epsilon\boldsymbol{I})^{\frac{p}{2}-1}\boldsymbol{V}^*. \tag{44}$$

One choice of the matrix square root $\sqrt{\boldsymbol{H}}$ is

$$\sqrt{\boldsymbol{H}} = \boldsymbol{V}(\boldsymbol{\Lambda}+\epsilon\boldsymbol{I})^{-\frac{q}{2}} = [(\lambda_1+\epsilon)^{-\frac{q}{2}}\boldsymbol{v}_1, \ldots, (\lambda_M+\epsilon)^{-\frac{q}{2}}\boldsymbol{v}_M]. \tag{45}$$

where $q = 1-\frac{p}{2}$. However, the least-squares subproblem of the algorithm only needs as input the filter $\boldsymbol{h}$ defined in (35), which is determined by the columns $\boldsymbol{h}_i$ of $\sqrt{\boldsymbol{H}}$. Making the substitutions $\boldsymbol{h}_i = (\lambda_1+\epsilon)^{-\frac{q}{2}}\boldsymbol{v}_i$ in (35) gives the following update:

$$\boldsymbol{h} = \sum_{i=1}^N (\lambda_i+\epsilon)^{-q}(\widetilde{\boldsymbol{v}}_i * \boldsymbol{v}_i); \quad q = 1-\frac{p}{2}. \tag{46}$$

Note that the weight $(\lambda_i+\epsilon)^{-q}$, $q>0$, is large only when the filter $\boldsymbol{v}_i$ is close to the null space of $\breve{\mathcal{T}}(\boldsymbol{x})$, i.e., when $\boldsymbol{v}_i$ is an annihilating filter. Therefore, $\boldsymbol{h}$ can be thought of as a weighted average of all the annihilating filters for $\breve{\mathcal{T}}(\boldsymbol{x})$. Notice that the only effect of changing the Schatten-$p$ penalty parameter $0 \le p \le 1$ is to change the exponent $q$ in the computation of (46). Smaller values of $p$ (larger values of $q$) are more likely to promote low-rank solutions at the iterations progress, since the weights on the filters close to the null space will be higher; see Figure 4 for an illustration of this effect.

---

**Algorithm 2:** GIRAF-$p$ algorithm for structured low-rank matrix recovery

---

Initialize $\boldsymbol{x}^{(0)}$ and choose $\epsilon_0 > 0$;
**for** $n = 1$ *to* $I_{max}$ **do**
- **Step 1: Annihilating Filter Update**
  - Build $\boldsymbol{G} = \breve{\mathcal{T}}(\boldsymbol{x}^{(n-1)})^*\breve{\mathcal{T}}(\boldsymbol{x}^{(n-1)})$ using (43);
  - Find eigendecomposition $(\lambda_i, \boldsymbol{v}_i)_{i=1}^N$ of $\boldsymbol{G}$;
  - Compute re-weighted annihilating filter
    $\boldsymbol{h} = \sum_{i=1}^N (\lambda_i + \epsilon_{n-1})^{\frac{p}{2}-1}(\widetilde{\boldsymbol{v}}_i * \boldsymbol{v}_i)$;
  - Convert filter to weights $\boldsymbol{D} = \operatorname{diag}(\boldsymbol{F}^*\boldsymbol{P}_{2\Lambda}^*\boldsymbol{h})$
- **Step 2: Least-squares Annihilation**
  - Solve least-squares problem:
  - $\boldsymbol{x}^{(n)} = \arg\min_{\boldsymbol{x}} \|\boldsymbol{A}\boldsymbol{x}-\boldsymbol{b}\|^2 + \lambda C_p \|\boldsymbol{D}^{\frac{1}{2}}\boldsymbol{F}^*\boldsymbol{M}\boldsymbol{x}\|^2$
  - by ADMM iterations (37)–(39) or CG;
- Choose $\epsilon_n$ such that $0 < \epsilon_n \le \epsilon_{n-1}$;

**end**

---

*D. Extension to liftings with multiple blocks*

To simplify the derivation of the GIRAF algorithm we assumed that the original matrix lifting $\mathcal{T}(\boldsymbol{x})$ consists of a single multi-level Toeplitz block. However, the GIRAF algorithm can easily be modified to accommodate liftings with a vertical block structure, as in (19), by applying the half-circulant approximation to each block. In this case, one can show the GIRAF least-squares problem simplifies to

$$\min_{\boldsymbol{x}} \|\boldsymbol{A}\boldsymbol{x}-\boldsymbol{b}\|_2^2 + \lambda C_p \sum_{i=1}^K \|\boldsymbol{D}^{\frac{1}{2}}\boldsymbol{F}^*\boldsymbol{M}_i\boldsymbol{x}\|_2^2. \tag{47}$$

where $\boldsymbol{D} = \operatorname{diag}(\boldsymbol{d})$, and the annihilation weights $\boldsymbol{d}$ are computed the same as in (33) and (46) from an eigendecomposition of the Gram matrix $\boldsymbol{G} = \sum_{i=1}^K \breve{\mathcal{T}}_i(\boldsymbol{x})^*\breve{\mathcal{T}}_i(\boldsymbol{x})$.

### E. Implementation Details

An overview of the GIRAF algorithm is given in Algorithm 2. Empirically we find that several of the heuristics in [45] for setting the smoothing parameter $\epsilon_n$ work well for the GIRAF algorithm, as well. In particular, we set $\epsilon_0 = \lambda_{\text{max}}/100$ where $\lambda_{\text{max}}$ is the maximum eigenvalue of the Gram matrix $\breve{\mathcal{T}}(\boldsymbol{x}^{(0)})^*\breve{\mathcal{T}}(\boldsymbol{x}^{(0)})$, which is obtained as a by-product of the first iterate of the algorithm. We recommend decreasing the smoothing parameter $\epsilon_n$ exponentially as $\epsilon_n = \epsilon_0(\eta)^{-n}$, where $\eta > 1$ is a fixed parameter. We find that $\eta = 1.2$ is suitable for a wide range of problem instances. The ADMM approach for solving the GIRAF least-squares subproblem also requires a conditioning parameter $\delta$. We typically set $\delta = 10$, which was found to represent a desirable trade-off between speed and accuracy for a wide range of problem instances; see Figure 7.

## V. Numerical experiments

In this section we perform several experiments to investigate properties of the GIRAF algorithm, and compare with other algorithms for structured low-rank matrix recovery. All experiments were conducted in MATLAB 8.5.0 (R2015a) on a Linux desktop computer with a Intel Xeon 3.20GHz CPU and 24 GB RAM.

### A. Behavior of GIRAF algorithm

*1) Half-circulant approximation, oversampled grid, and choice of $p$:* The GIRAF algorithm relies on a half-circulant approximation (26) of Toeplitz matrix liftings. Here we investigate the error induced by this approximation in a 1-D signal recovery setting, where the matrix lifting is single-level Toeplitz: $\mathcal{T}(\boldsymbol{x}) = \mathsf{Toep}(\boldsymbol{x})$. We randomly undersample the exact Fourier coefficients $\boldsymbol{x}_0 = (\hat{\rho}[k]\colon |k| < 64)$ of a 1-D stream of $r$ Diracs $\rho(x)$ as in (1), and attempt to recover the missing samples by applying the IRLS-$p$ to the Toeplitz matrix lifting $\mathcal{T}(\boldsymbol{x})$ having dimensions $113 \times 15$, and by using the GIRAF-$p$ algorithm, which is equivalent to applying IRLS-$p$ to the half-circulant approximated lifting $\breve{\mathcal{T}}(\boldsymbol{x})$ having dimensions $157 \times 15$. We observe the IRLS-0 algorithm generally outperforms the IRLS-1, consistent with the results in [21]. Furthermore, the quality of the GIRAF-$p$ reconstructions mimics those obtained with IRLS-$p$, but with small approximation errors. However, the approximation error is noticeably less severe in the non-convex $p = 0$ case.

We investigate this phenomenon in another experiment shown in Figure 6. Using the same experimental setup, we varied the working grid size of the problem, i.e., we varied the size of the oversampled reconstruction grid $\Delta'$ defined in Section IV-A; this only changes the number of rows $M$ in the lifting $\breve{\mathcal{T}}(\boldsymbol{x})$. We vary $M = 127, ..., 255$, corresponding to a maximum oversampling factor of 2. For the convex GIRAF-1 algorithm, we find that the reconstruction error stagnates with respect to the oversampling factor. However, for the non-convex GIRAF-$p$ algorithms, $p < 1$, the reconstruction diminishes significantly with the oversampling factor. We additionally perform a similar experiment using the gradient weighted lifting scheme (9) on synthetic data that is piecewise constant and known to be low-rank in the lifted domain (the **SL** dataset shown in Figure 9) with data size $201 \times 201$, filter size $25 \times 25$, and $50\%$ random samples. In this case we also obtain similar results—the nonconvex GIRAF algorithms achieve substantially smaller approximation errors that decrease with the oversampling factor. While this suggests one should use a very large oversampled grid to eliminate approximation errors, there is a trade-off with computational cost, since the oversampled grid represents the working dimension of the problem. The extent of oversampling will depend on the problem setting and required precision.

These experiments suggest it is better to use a non-convex GIRAF-$p$, $0 \leq p < 1$, algorithm over the convex GIRAF-1 algorithm in order to mitigate errors due to the half-circulant approximation.

*2) ADMM approach to GIRAF least-squares solution:* In Figure 7 we compare the computation time of the proposed ADMM approach (37)–(39) for solving the GIRAF least-square subproblem (34) against standard CG and LSQR solvers. Here we investigate the GIRAF-0 algorithm applied to a Fourier domain recovery experiment using the gradient weighted lifting scheme (9). We take the measurement operator $\boldsymbol{A}$ to be a Fourier domain sampling, and sample 25% of the data uniformly at random, using the Shepp-Logan phantom (data size $256 \times 256$ and filter size of $35 \times 35$). The metric used to evaluate each approach is the normalized mean-square difference NMSD $= \|\boldsymbol{x} - \boldsymbol{x}_*\|^2 / \|\boldsymbol{x}_*\|^2$ where $\boldsymbol{x}$ is the current iterate, and $\boldsymbol{x}_*$ is the true solution, which was obtained by running CG algorithm to high precision. The proposed ADMM approach shows nearly a ten-fold increase in over CG and LSQR. The rate at which the ADMM approach reduces the NMSD shows some sensitivity to the parameter $\delta$, but only in the high-accuracy regime (NMSD $< 10^{-4}$).

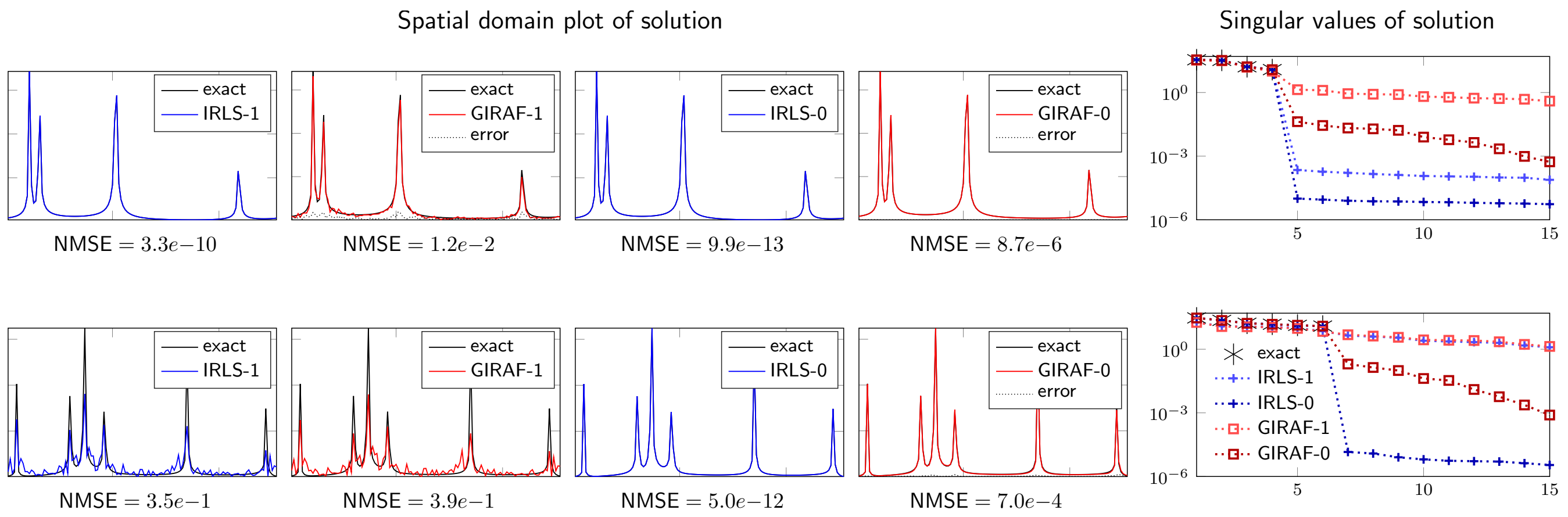

Fig. 5. Comparison of solutions obtained by the IRLS-$p$ and the GIRAF-$p$ algorithm, $p = 1, 0$, for recovery of Fourier coefficients of $r$ Diracs, as in (1), from non-uniform random Fourier samples. The top row shows an experiment ($r = 4$, 50% random samples) where both IRLS-1 and IRLS-0 succeed in recovering the signal. In this case the GIRAF-1 algorithm shows non-negligible approximation errors, but the GIRAF-0 result is close to exact. The bottom row shows an experiment ($r = 6$, 33% random samples) where IRLS-1 fails, but IRLS-0 is successful; similarly, the GIRAF-1 recovery fails, while the GIRAF-0 result is again close to exact. Shown in the right column are the singular values of the exact matrix lifting $\mathcal{T}(\boldsymbol{x}^\star)$ where $\boldsymbol{x}^\star$ is the solution obtained from each algorithm. Note GIRAF-0 solutions are still approximately low-rank under the exact lifting.

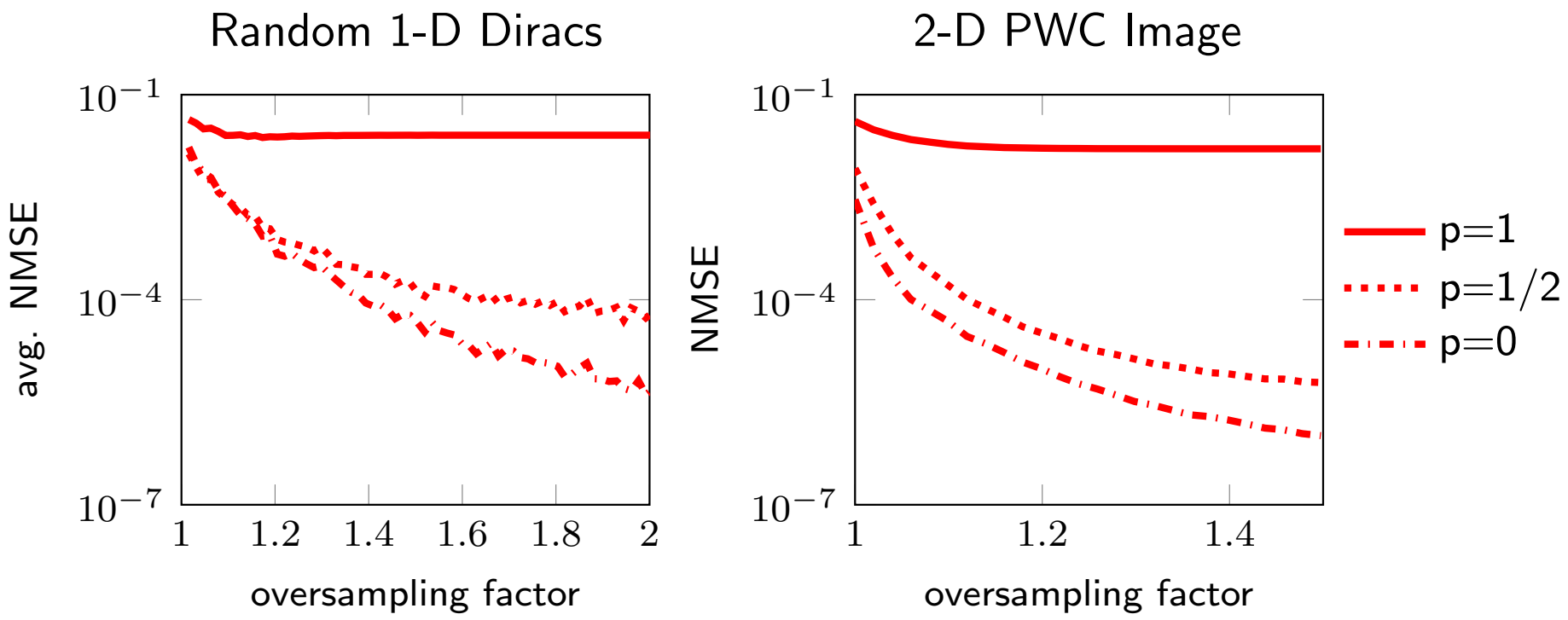

Fig. 6. Reconstruction error versus oversampling factor for GIRAF-$p$ algorithm on the recovery of synthetic low-rank data. The left plot corresponds to the recovery Fourier coefficients of Dirac in 1-D, using the same settings as in Figure 5; we report the average NMSE over 50 random trials. The right plot shows recovery of Fourier coefficients of a synthetic 2-D piecewise constant image using the gradient weighted lifting (9). Observe that the NMSE stagnates for the $p = 1$ case (convex penalty). However, the NMSE diminishes towards zero with the oversampling factor when $p = 0.5$ and $p = 0$ (non-convex penalties). This suggests it is important to use the non-convex versions of GIRAF on a sufficiently large oversampled grid to mitigate errors due to the half-circulant approximation.

### B. Recovery Experiments

To demonstrate the benefits of the GIRAF algorithm for structured low-rank matrix recovery, we focus on the problem of undersampled MRI reconstruction in 2-D. In this setting, the goal is to recover an array $\boldsymbol{x}_0$ of the Fourier coefficients of an image from missing or corrupted Fourier samples. In addition to the direct IRLS-$p$ algorithm (see Algorithm 1), we compare GIRAF against the following algorithms proposed for structured low-rank matrix recovery:

- **Alternating projections (AP)**. Also known as Cadzow's method [58], this approach was adopted for noisy finite-rate-of-innovation (FRI) signal recovery in [38], [59], [60], and for auto-calibrated parallel MRI reconstruction in [2]. Reconstruction in an AP algorithm is posed as

$$\min_{\boldsymbol{x}} \|\boldsymbol{A}\boldsymbol{x} - \boldsymbol{b}\|^2 \quad \text{subject to} \quad \begin{cases} \boldsymbol{X} = \mathcal{T}(\boldsymbol{x}) \\ \operatorname{rank} \boldsymbol{X} \leq r, \end{cases} \tag{48}$$

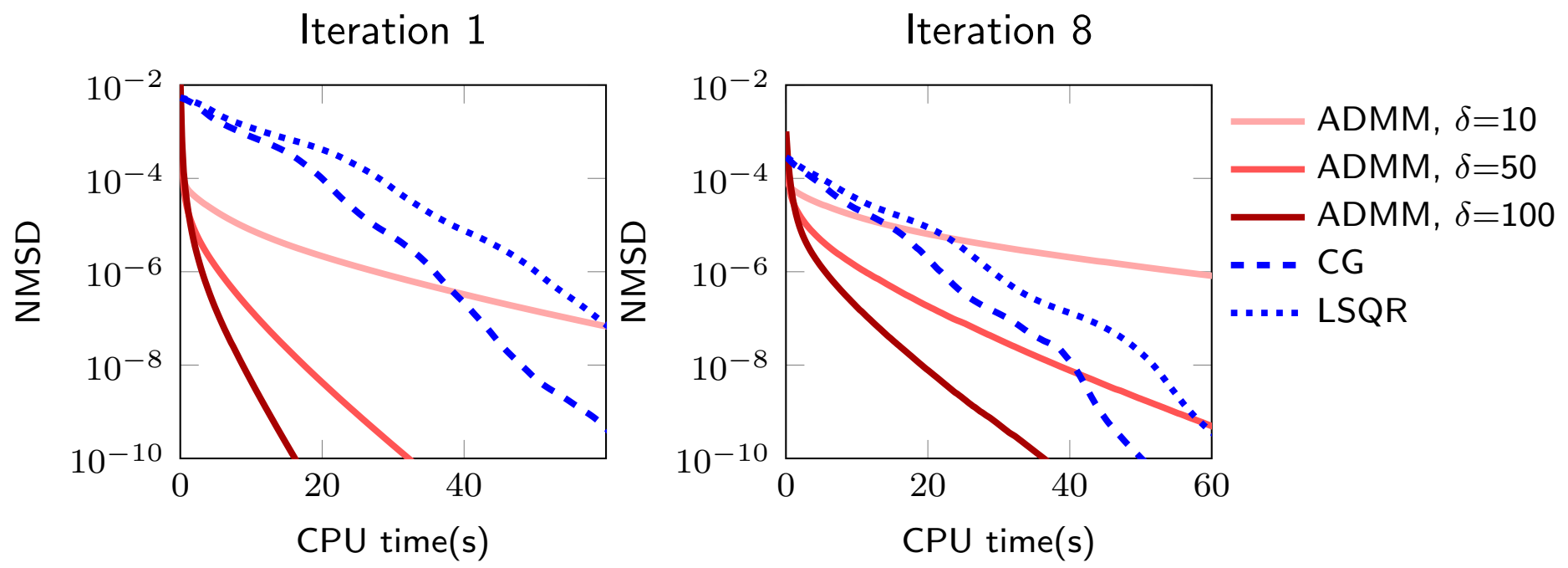

Fig. 7. Computation time comparison for solvers of the GIRAF least-squares problem. We study the least-squares problem arising in iterations 1 and 8 of the GIRAF algorithm applied to a Fourier recovery experiment. The proposed ADMM approach shows nearly ten-fold increase in the rate at which NMSD is decreased over CG and LSQR, provided the pre-conditioning parameter $\delta$ is chosen correctly.

  which is solved by alternately projecting $\boldsymbol{X}$ onto: (1) the set of matrices with rank less than or equal to $r$, (2) the space of linear structured matrices specified by the range of $\mathcal{T}$, and (3) the data fidelity constraint set. Note that the AP algorithm requires a rank parameter $r$. A novel extension of the AP algorithm incorporating proximal-smoothing is used in the LORAKS framework [4], [10], [61] for structured low-rank based MRI reconstruction. We use **AP-PROX** (AP with proximal smoothing) to refer to the algorithm introduced in [4], to distinguish it from the structured low-rank matrix models also introduced in [4]; see the Supplementary Materials for more details.
- **Singular value thresholding (SVT)**. Proposed in [62] for general low-rank matrix recovery by nuclear norm minimization, and adapted to Hankel structured matrix case in [45], and for 2-D spectral estimation in [3]. The SVT algorithm minimizes the objective (12) or (16) (with $p = 1$) via iterative terative soft thresholding of the singular values of the lifted matrix.
- **Singular value thresholding with factorization heuristic (SVT+UV)**. This approach was proposed in [63] for general structured low-rank matrix recovery, and was adopted by the present authors in [8] for structured low-rank based super-resolution MRI reconstruction. It is also adopted in the ALOHA framework [31] for a variety of imaging applications, including structured low-rank based MRI reconstruction. The SVT+UV approach also minimizes the objective (12) or (16) (with $p = 1$), but uses the matrix factorization characterization of the nuclear norm:

$$\|\boldsymbol{X}\|_* = \min_{\substack{\boldsymbol{U} \in \mathbb{C}^{M \times R} \\ \boldsymbol{V} \in \mathbb{C}^{N \times R} \\ \boldsymbol{X} = \boldsymbol{U}\boldsymbol{V}^*}} \frac{1}{2} \left( \|\boldsymbol{U}\|_F^2 + \|\boldsymbol{V}\|_F^2 \right), \tag{49}$$

  which holds true provided the inner dimension $R$ of $\boldsymbol{U}$ and $\boldsymbol{V}^*$ satisfies $R \geq \operatorname{rank} \boldsymbol{X}$. This allows efficient solutions of the SVT subproblems by matrix inverses; see the Supplementary Materials for more details.

To aid in reproducibility of our results, we give the implementation details and pseudo-code for all these algorithms in the Supplementary Materials. We note only the AP, AP-PROX, and SVT+UV algorithms were previously used for the type of large-scale SLRA problems that we consider. However, we also include comparisons against SVT and IRLS for benchmark purposes.

To compare the recovery performance of each algorithm, we measure the quality of each reconstruction by its normalized mean square error $\text{NMSE} = \|\boldsymbol{x} - \boldsymbol{x}_0\|^2/\|\boldsymbol{x}_0\|^2$, where $\boldsymbol{x}$ is the solution obtained from the algorithm, and $\boldsymbol{x}_0$ is the ground truth data. We use the NMSE as a error metric since the competing algorithms minimize different cost functions, and cannot be compared on the basis of how well they minimize a common cost function. NMSE is also a commonly used error metric[4] in the MRI reconstruction literature (see, e.g., [1]).

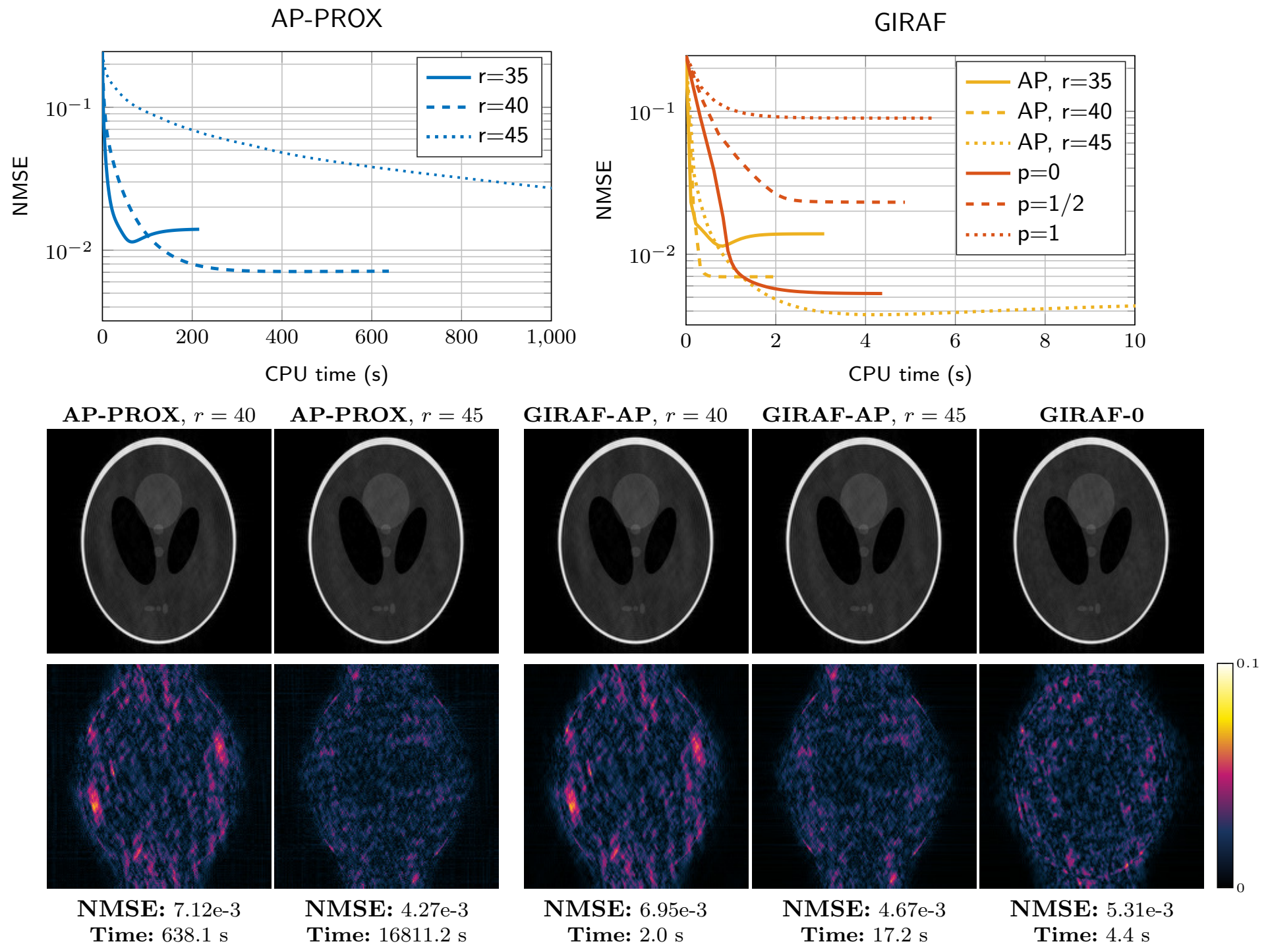

Fig. 8. Comparison of GIRAF with AP-PROX algorithm for recovery with the C-LORAKS SLRA model. Plotted is the per iteration NMSE against elapsed CPU time in the recovery of synthetic data from the indicated from a uniform random Fourier sampling pattern.

*1) GIRAF for LORAKS recovery:* First, to demonstrate the benefit of the GIRAF algorithm for an existing SLRA approach in MRI, we apply GIRAF to the LORAKS constrained MRI reconstruction framework [4]. In Figure 8 we compare the performance of the GIRAF to the AP-PROX algorithm proposed in [4], [64], using publicly available code and data[5]. For simplicity we restrict our comparisons to the "C-LORAKS" matrix lifting, which is based on a spatial sparsity assumption; this is equivalent to a single-block two-level Toeplitz lifting (24) with $\boldsymbol{M}$ the identity matrix. The LORAKS framework assumes filters with approximately circular support and we implement filters with the same size and support in the GIRAF algorithm to ensure a fair comparison. For the AP-PROX algorithm we used the default parameters distributed with the LORAKS code. We tested both algorithms on the dataset provided with the LORAKS code ($180 \times 180$ Fourier coefficients of the Shepp-Logan phantom with phase) and one of the provided Fourier domain sampling masks (a uniform random sampling pattern accounting for $\approx 63\%$ of the data). We find the GIRAF algorithm converges an order of magnitude faster than the AP-PROX algorithm (2-5 s versus 50-200+ s) and, in these experiments, to a solution with lower relative error as measured by the NMSE. In addition to testing the default LORAKS rank cutoff parameter $r = 35$, we also varied $r$ to study the effect on the reconstruction. We observe that the AP-PROX algorithm shows significant variation in final NMSE depending on the rank estimate $r$. Similar behavior was observed in the experiments in [5].

[4]In place of NMSE, some authors use the signal-to-noise ratio (SNR) instead, defined as $-10\log_{10}(\text{NMSE})$.
[5]http://mr.usc.edu/download/loraks/

| Problem | | | | | AP | | SVT | | SVT+UV | | IRLS-0 | | GIRAF-0 | |
|---|---|---|---|---|---|---|---|---|---|---|---|---|---|---|
| dataset | data size | filter | rank | USF | # iter | time (s) | # iter | time (s) | # iter | time (s) | # iter | time (s) | # iter | time (s) |
| **PWC1** | 65×65 | 9×9 | 32 | 0.50 | 17 | 5.4 | 11 | 1.1 | 11 | 0.8 | 3 | 17.5 | 3 | **0.6** |
| **PWC1** | 65×65 | 9×9 | 32 | 0.33 | 39 | 12.8 | - | Inf | 19 | **1.4** | 4 | 24.5 | 5 | **1.4** |
| **PWC2** | 129×129 | 17×17 | 118 | 0.50 | 19 | 48.0 | 9 | 25.2 | 15 | 21.6 | 4 | 566.6 | 3 | **1.1** |
| **PWC2** | 129×129 | 17×17 | 118 | 0.33 | 90 | 228.0 | - | Inf | 24 | 36.2 | 7 | 1032.8 | 5 | **2.6** |
| **SL** | 201×201 | 25×25 | 336 | 0.65 | 31 | 762.0 | 13 | 338.7 | 11 | 150.9 | - | Mem | 3 | **3.7** |
| **SL** | 201×201 | 25×25 | 336 | 0.50 | 88 | 2151.7 | - | Inf | 36 | 494.4 | - | Mem | 4 | **4.7** |
| **BRAIN** | 255×255 | 45×45 | 1250 | 0.65 | - | Mem | - | Mem | 23 | 4270.3 | - | Mem | 3 | **25.2** |
| **BRAIN** | 255×255 | 45×45 | 1250 | 0.50 | - | Mem | - | Mem | - | Inf | - | Mem | 5 | **55.3** |

**Stopping criteria**: NMSE $\leq 10^{-4}$. **Key**: Inf - Algorithm converged above NMSE threshold. Mem - Not enough memory to run algorithm. **Computer specs**: Intel Xeon 3.20GHz CPU and 24 GB RAM.

TABLE I
COMPUTATION TIME OF ALGORITHMS FOR SLRA RECOVERY OF PIECEWISE CONSTANT IMAGES.

To investigate whether the improvement offered by GIRAF in this setting is due to a difference in the cost functions between the GIRAF-$p$ approach and AP-PROX algorithm, we also adapt the GIRAF approach to to minimize a cost-function similar to AP-PROX. This was done by redefining the matrix $\boldsymbol{H}$ in (44) to be a projection onto the eigenspace of the Gram matrix $\boldsymbol{G} = \breve{\mathcal{T}}(\boldsymbol{x})^*\breve{\mathcal{T}}(\boldsymbol{x})$ determined by its $r$ smallest eigenvalues. Namely, we set $\boldsymbol{H} = \boldsymbol{V}\boldsymbol{Q}_r\boldsymbol{V}^*$, where $\boldsymbol{V}$ is a basis of eigenvectors of $\boldsymbol{G}$, and $\boldsymbol{Q}_r$ is a diagonal matrix with ones along the diagonal at the locations of the $r$ smallest eigenvalues, and zeros elsewhere. We label this approach GIRAF-AP, and compare GIRAF-AP with AP-PROX at various values of rank cutoffs. We find that GIRAF-AP and AP-PROX give similar results in terms of NMSE for the different values of $r$.

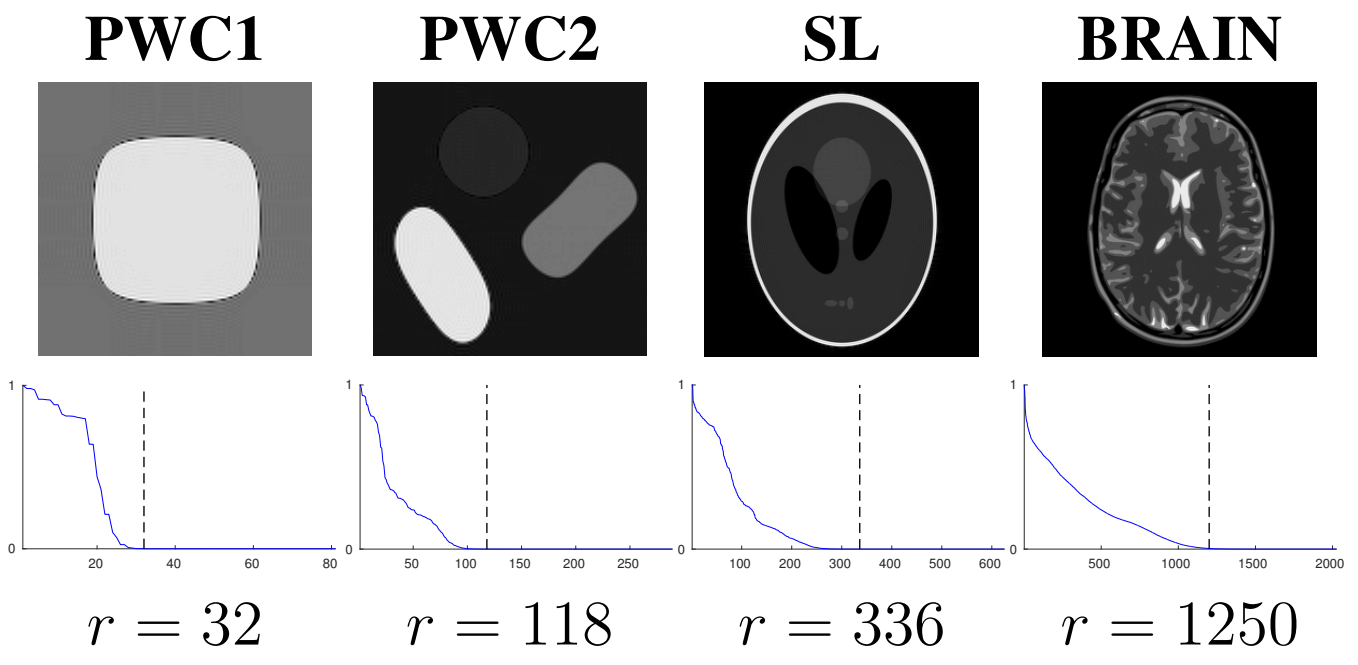

Fig. 9. Synthetic piecewise constant images used in structured low-rank recovery experiments. Shown below the images are the normalized singular value spectrum of the gradient weighted lifted matrix constructed from the fully sampled Fourier coefficients, where the dotted line indicates the rank estimate $r$.

*2) GIRAF for piecewise constant image recovery:* Here we test the GIRAF algorithm for the recovery of images belonging to the piecewise constant SLRA model proposed in [9], which uses the gradient weighted matrix lifting described in (9). We adapt the AP, AP-PROX, SVT, SVT+UV, and IRLS algorithms to this setting. To compare computation time among algorithms, we first experiment with recovering simulated piecewise constant images from their undersampled Fourier coefficients. The datasets used in these experiments are shown in Figure 9. In each experiment, we sample uniformly at random in Fourier domain at the specified undersampling factor (USF). Since the simulated Fourier coefficients are noise-free, each algorithm incorporates equality data constraints, as in the formulation (12). To ensure a fair comparison, the algorithms that use a rank estimate (AP, SVT+UV) were passed either the exact rank of the simulated data, which was calculated either using the known bandwidth of the level-set polynomial describing the edge-set of the image (see [9]), or the index for which the normalized singular values were less than $10^{-2}$. The results of these experiments are shown in Table I. We report the CPU time and number of iterations for each algorithm to reach NMSE $\leq 10^{-4}$. Observe that for small to medium problem sizes (**PWC1**, and **PWC2**), the GIRAF algorithm is competitive or significantly faster than state-of-the-art methods. For the large-scale problems (**SL**, and **BRAIN**) the GIRAF algorithm converges orders of magnitude faster than competing algorithms, demonstrating its superior scalability. GIRAF is also successful on all the "hard" problem instances where SVT fails to converge below the set NMSE tolerance.

In Figure 10 we show the results of a similar recovery experiment, but where the Fourier samples are corrupted with noise. We test on the **SL** dataset with USF $= 0.65$, adding complex white Gaussian noise to the Fourier samples such that the signal-to-noise ratio (SNR) is approximately 22 dB. Here we test GIRAF against AP-PROX and SVT-UV in their regularized formulations, similar to (16), and we tune regularization and rank parameters to obtain the optimal NMSE in each case. We use a filter of size $21 \times 21$ in all algorithms. In Figure 10 we report the time each algorithm took to reach $1\%$ of the final NMSE, where the final NMSE is obtained by running each algorithm until the relative mean squared difference between iterates is less than $10^{-8}$. We observe the GIRAF-0 algorithm shows similar runtime as in the noise-free setting, and converges to a solution with smaller NMSE. GIRAF-0 also converged to a solution to smaller NMSE compared with GIRAF-1/2 and GIRAF-1.

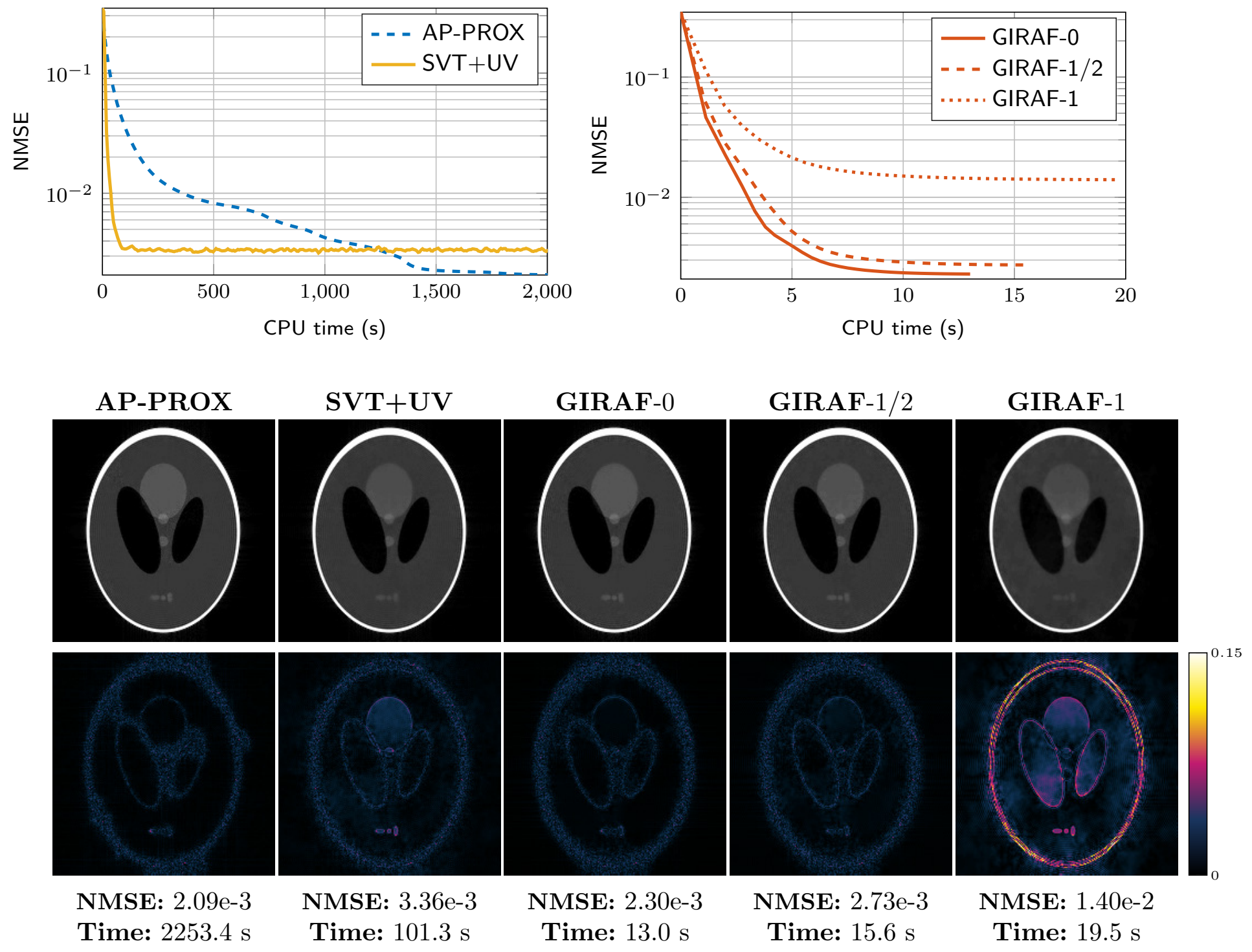

Fig. 10. Comparison of GIRAF with competing algorithms for structured low-rank matrix recovery with noisy data using a piecewise constant SLRA model. Plotted is the per iteration NMSE against elapsed CPU time in the recovery of synthetic data from noisy random uniform Fourier samples (USF=0.65, sample SNR=22dB).

*3) Application to compressed sensing MRI reconstruction:* We demonstrate the GIRAF algorithm for the recovery of real MRI data from undersampled Fourier data using an SLRA approach. The datasets we experiment on were obtained from a fully sampled four-coil parallel MRI acquisition of a human brain. We compressed the data into a single virtual coil using an SVD-based technique [65], then retrospectively undersampled the virtual coil data. Each aquisition consists of $256 \times 170$ pixels. To simulate a "calibrationless" sampling strategy, similar to those investigated in [2], [4], we sample in Fourier domain uniformly at random.

In Figure 11(a) we show the results of recovering undersampled data using the C-LORAKS SLRA model with GIRAF, SVT+UV, and AP-PROX. To demonstrate the potential benefit of a Fourier domain SLRA approach over standard discrete spatial domain penalties, we also compare against a total variation (TV) regularized reconstruction implemented with an efficient ADMM/Split-Bregman algorithm [66]. We use a filter of size $11 \times 11$, and for each algorithm, we tuned the regularization and rank parameters to maximize SNR $= -10\log_{10}(\text{NMSE})$. We ran each algorithm until the relative mean squared difference between iterates was less than $10^{-8}$. On this dataset the GIRAF algorithm converged in 1-3 s compared to over a minute using SVT+UV or AP-PROX. Here SVT+UV gives slightly better SNR than GIRAF-0 (+0.3dB), but visually they appear indistinguishable. The GIRAF-1/2 and GIRAF-1 reconstructions compare poorly to GIRAF-0, SVT+UV, AP-PROX, similar to our experiments on synthetic data. The TV reconstruction takes under a second, but gave poor quality results in this setting, similar to the results obtained in [4].

In Figure 11(b) we show the results of recovering undersampled data using the piecewise constant SLRA model. Here the filter size was set to $21 \times 21$; the larger filter size is necessary to obtain optimal results with the piecewise constant model (see Figure 12). The GIRAF reconstruction took 7-10 s, while SVT+UV and AP-PROX took several minutes. The GIRAF-0 reconstruction matches SVT+UV in SNR and visual quality, but is computed orders of magnitude faster. Similarly, AP-PROX and GIRAF-1/2 give comparable results, while GIRAF-1 performs poorly. [6] Observe that the best SLRA reconstructions show more than 3dB improvement in SNR over the TV reconstruction, demonstrating clear benefits of an SLRA approach in this setting. However, we note that we use a relatively high undersampling factor to illustrate the difference in reconstruction quality between competing methods; the reconstruction quality may not be appropriate for radiological evaluations.

One advantage of the GIRAF algorithm is that it scales well with filter size over competing algorithms. We show in Figure 12 that using larger filter sizes directly translates to improved image quality in the case of the gradient weighted lifting with modest increases in runtime. Here we re-ran the experiment in Fig. 11(b) but with different filter-sizes, which shows that substantial improvement in SNR can be from increasing the filter size from $11\times 11$ to $21\times 21$. Additional improvement in SNR is attained by increasing the filter to $31 \times 31$, but at the expense of a longer reconstruction time.

## VI. Discussion and conclusion

We introduced the GIRAF algorithm for multi-level Toeplitz/Hankel low-rank matrix recovery problems arising in MRI reconstruction and other imaging applications. The algorithm is based on the IRLS approach for low-rank matrix completion, combined with a novel half-circulant approximation of multi-level Toeplitz matrices. This approximation dramatically reduces the computational complexity and memory demands of the direct IRLS algorithm for large problem sizes. Due to the monotonicity of the Schatten-$p$ regularizer with respect to the half-circulant approximation, we show it suffices to regularize the approximated matrix, rather than the original multi-level Toeplitz matrix. Compared to previous approaches the GIRAF algorithm is an order of magnitude faster on several realistic problem instances. This is because in the GIRAF algorithm performs most of its computations in the original "un-lifted" problem domain with FFTs, rather than in the "lifted" matrix domain. An important feature of the GIRAF algorithm is that it can accommodate larger filters required by more sophisticated image priors, such as the off-the-grid piecewise constant image model proposed in [9]. This enables SLRA/annihilating filter approaches for a much wider class of imaging problems, including realistic large-scale multidimensional datasets encountered in MRI reconstruction tasks.

Unlike some previous approaches, the GIRAF algorithm does not require strict low-rank approximations, nor does it require an estimate of the underlying rank or model-order. Arguably the optimal selection of the $\epsilon_0$ smoothing parameter used in the GIRAF algorithm implicitly performs model-order selection. However, we propose a method to automatically set this parameter based on spectral properties of the initialization to the algorithm. Similar strategies could be used for other algorithms that require a rank estimate to reduce their dependence on the rank.

[6] The worse performance AP-PROX in this setting may be due to slow convergence, which limited our ability to tune rank and regularization parameters.

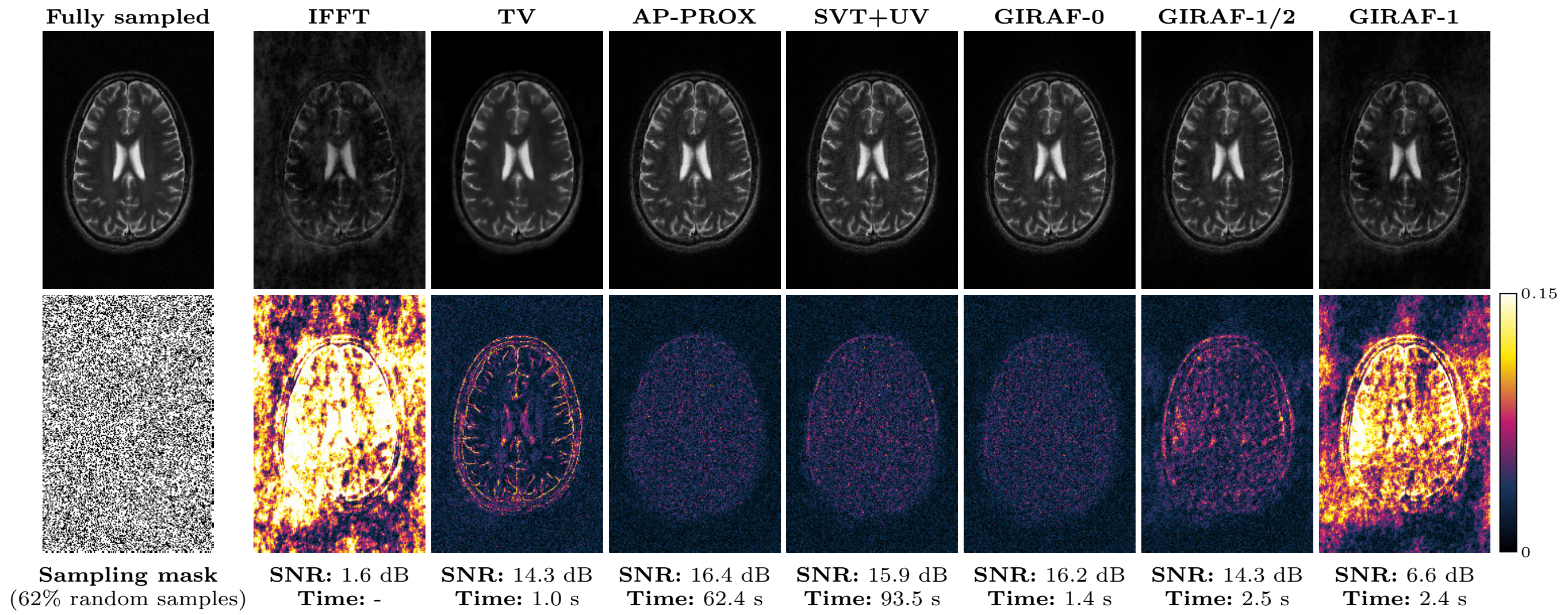

(a) C-LORAKS SLRA model

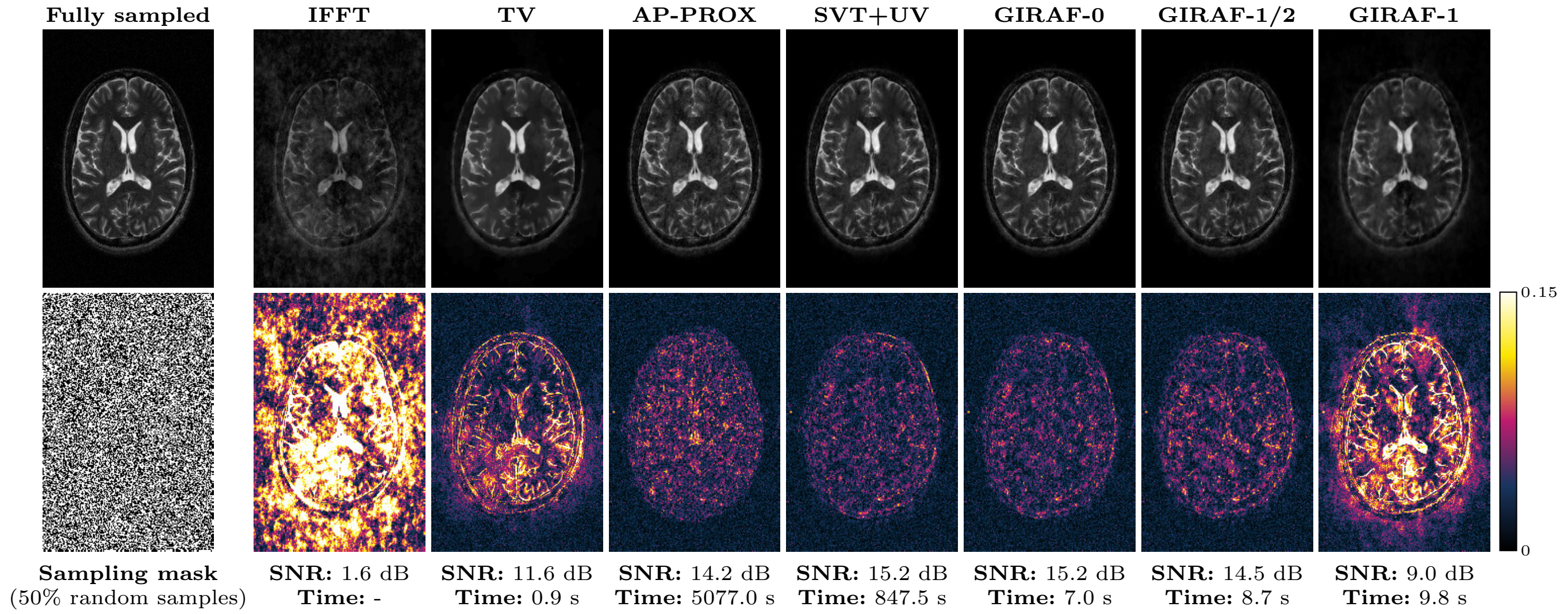

(b) Piecewise constant SLRA model

Fig. 11. Undersampled calibrationless MRI reconstruction on real data using SLRA models: (a) using a C-LORAKS SLRA model with filter size $11 \times 11$ with USF=0.62, and (b) using a piecewise constant SLRA model with a filter size $21 \times 21$ USF=0.50. For SLRA recovery we compare algorithms AP-PROX, SVT+UV, and GIRAF-$p$ with $p = 0, 1/2, 1$. We also include a zero-filled IFFT of the undersampled data and a total variation (TV) regularized reconstruction for comparison. Error images are shown in the bottom row.

The multi-level Toeplitz matrix lifting model considered in this work may seem narrow in application, but in fact a variety of existing SLRA models belong to this class, giving the GIRAF algorithm fairly wide applicability. Single-level Toeplitz/Hankel matrix liftings form the foundation for several modern approaches in spectral estimation [60], direction-of-arrival estimation [67], and system identification [45]. Multi-level Toeplitz/Hankel matrix liftings have been used for multi-dimensional spectral estimation [3], [68], and video inpainting [31], [69]. Of particular interest is the application of the framework to super-resolution microscopy [70], which is of high significance in biological imaging. SLRA methods have been seen to be very promising in this context [30], offering more accurate localizations than state-of-the-art methods. The GIRAF algorithm could potentially enable the extension of these methods to recover the whole image rather than small patches, due to its low memory demand. Similarly, most of the SLRA models proposed for MRI reconstruction have a block multi-level Toeplitz/Hankel matrix structure. For example, the SAKE framework for parallel MRI reconstruction [2] considers a matrix lifting with a horizontal block structure, where each block is a multi-level Hankel matrix. Similarly, the matrix liftings considered in the LORAKS [4], [10] and ALOHA [6] frameworks have block multi-level Toeplitz/Hankel structure. In this work we restricted our attention to matrix liftings having a vertical block structure. While our current algorithm is not readily applicable to these settings, the GIRAF framework could also potentially be generalized to more general block-rectangular liftings, such as those used parallel MRI context [2], [10]; we save this as a topic for future work.

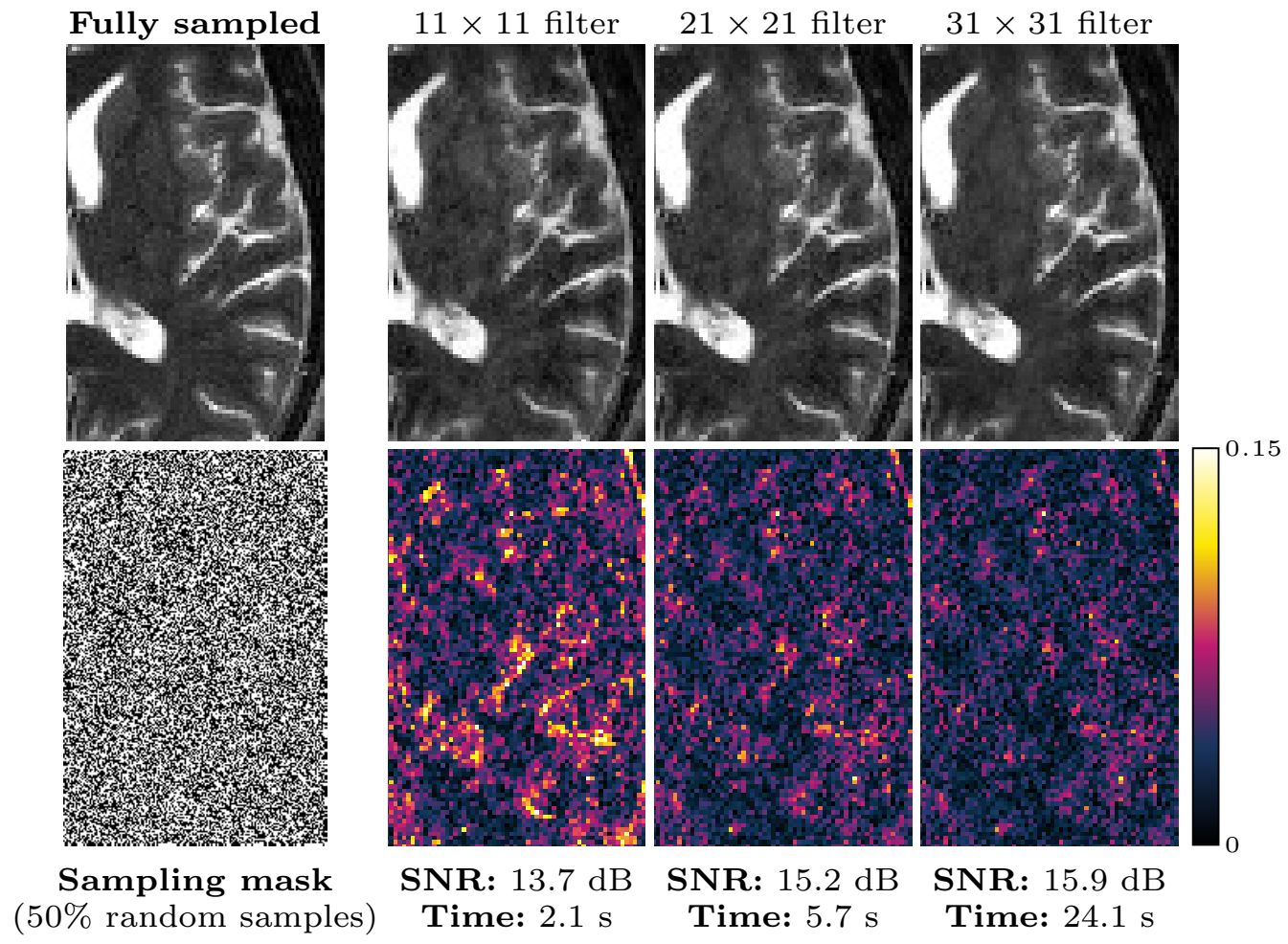

Fig. 12. Effect of filter size on reconstruction; zoomed for detail. We reconstruct the data shown in Fig. 11(b) using the GIRAF-0 algorithm and the piecewise constant SLRA model. The reconstruction shows a 2 dB improvement in SNR by increasing the filter size from $11 \times 11$ to $31 \times 31$, with only a modest increase in runtime (2 s versus 24 s).

Finally, the approximation intrinsic to the GIRAF approach means that the algorithm might not be appropriate for all SLRA problems. For instance, if the rank of the lifted matrix is expected to be very small and is known in advance, singular value thresholding or alternating projection algorithms may be more efficient or more accurate. However, for the large-scale SLRA problems encountered in imaging, where the rank of the lifted matrix is typically unknown and possibly large, our numerical experiments show the GIRAF approach represents a desirable trade-off between accuracy and runtime versus other state-of-the-art algorithms.

# A Fast Algorithm for Convolutional Structured Low-Rank Matrix Recovery: Supplementary Materials

Greg Ongie*, *Student Member, IEEE*, Mathews Jacob, *Senior Member, IEEE*

## I. MAJORIZATION-MINIMIZATION FORMULATION OF IRLS-$p$ ALGORITHM

Majorization-minimization (MM) algorithms iteratively minimize a given cost function by minimizing a surrogate cost function that upper bounds the original at each iteration [1]. More precisely, if $f(\boldsymbol{x})$ is any real-valued cost function on a domain $\Omega$, a function $g(\boldsymbol{x};\boldsymbol{x}_0)$ defined on $\Omega\times\Omega$ is said to *majorize* $f(\boldsymbol{x})$ at the fixed value $\boldsymbol{x}_0\in\Omega$ if

$$g(\boldsymbol{x};\boldsymbol{x}_0)\geq f(\boldsymbol{x}),\quad \text{for all}\quad \boldsymbol{x}\in\Omega \tag{1}$$

$$g(\boldsymbol{x}_0;\boldsymbol{x}_0)= f(\boldsymbol{x}_0) \tag{2}$$

If (1) and (2) hold for all $\boldsymbol{x}_0\in\Omega$ then $g$ is called a *majorizer* for $f$. An MM algorithm finds a minimizer for $f$ by sequentially solving for the iterates

$$\boldsymbol{x}^{(n+1)}=\arg\min_{\boldsymbol{x}\in\Omega} g(\boldsymbol{x};\boldsymbol{x}^{(n)}). \tag{3}$$

Using properties (1) and (2) it is easy to show the cost function $f$ must monotonically decrease at each iteration. Furthermore, if $f$ is strongly convex, the iterates $\boldsymbol{x}^{(n)}$ are guaranteed to converge to the unique global minimizer of $f$; when convexity is violated, the iterates are still guaranteed to converge to a stationary point of $f$ [2].

We show that Algorithm 1 can be interpreted as an MM algorithm for the cost function

$$\min_{\boldsymbol{x}} \|\boldsymbol{A}\boldsymbol{x}-\boldsymbol{b}\|_2^2+\lambda\|\mathcal{T}(\boldsymbol{x})\|_{p,\epsilon}^p \tag{4}$$

where $\|\cdot\|_{p,\epsilon}^p$ *smoothed* Schatten-$p$ penalty defined as

$$\|\boldsymbol{X}\|_{p,\epsilon}^p := Tr[(\boldsymbol{X}^*\boldsymbol{X}+\epsilon\boldsymbol{I})^{p/2}]=\sum_{i=1}^{N}(\sigma_i(\boldsymbol{X})^2+\epsilon)^{p/2} \tag{5}$$

for $0<p\leq 1$, and

$$\|\boldsymbol{X}\|_{0,\epsilon}^0 := \frac{1}{2}\log\det(\boldsymbol{X}^*\boldsymbol{X}+\epsilon\boldsymbol{I})=\sum_{i=1}^{N}\frac{1}{2}\log(\sigma_i(\boldsymbol{X})^2+\epsilon) \tag{6}$$

where $\sigma_i(\boldsymbol{X})$ are the singular values of $\boldsymbol{X}$. Specifically, we construct majorizers for the smoothed Schatten-$p$ penalties. Let $H_+^N$ and $H_{++}^N$ be the set of complex $N\times N$ matrices that are positive definite and positive semi-definite, respectively. The majorizers are derived from the following trace inequalities:

G. Ongie is with the Department of Electrical Engineering and Computer Science, University of Michigan, Ann Arbor, MI, 48104 USA, and M. Jacob is with the Department of Electrical and Computer Engineering, University of Iowa, Iowa City, IA, 52245 USA (e-mail: gongie@umich.edu; mjacob@uiowa.edu)

This work is supported by grants NIH 1R01EB019961-01A1, NSF CCF-1116067, ONR N00014-13-1-0202, and ACS RSG-11-267-01-CCE.

**Proposition 1.** *If* $\boldsymbol{Y}_0 \in H_{++}^N$ *then for all* $\boldsymbol{Y} \in H_+^N$ *we have*

$$Tr[\boldsymbol{Y}^q] \leq Tr[\boldsymbol{Y}_0^q + q\boldsymbol{Y}_0^{q-1}(\boldsymbol{Y} - \boldsymbol{Y}_0)] \tag{7}$$

*for all* $q \in (0, 1]$*. If additionally* $\boldsymbol{Y} \in H_{++}^N$*, then*

$$Tr[\log(\boldsymbol{Y})] \leq Tr[\log(\boldsymbol{Y}_0) + (\boldsymbol{Y}_0)^{-1}(\boldsymbol{Y} - \boldsymbol{Y}_0)]. \tag{8}$$

*Proof.* These are special cases of Klein's inequality (see, e.g., Proposition 2.5.2 in [3]), which states that for any concave function $f : [0, \infty) \to \mathbb{R}$ differentiable on $(0, \infty)$, and any $\boldsymbol{Y} \in H_+^N$, $\boldsymbol{Y}_0 \in H_{++}^N$ we have

$$Tr[f(\boldsymbol{Y}) - f(\boldsymbol{Y}_0) - f'(\boldsymbol{Y}_0)(\boldsymbol{Y} - \boldsymbol{Y}_0)] \leq 0 \tag{9}$$

where $f(\boldsymbol{X}) := \sum_i f(\lambda_i)\boldsymbol{P}_i$ when $\boldsymbol{X}$ has the spectral decomposition $\boldsymbol{X} = \sum_i \lambda_i \boldsymbol{P}_i$. Choosing $f(t) = t^q$ and $f(t) = \log(t)$ gives the desired inequalities. □

Let $\boldsymbol{X}, \boldsymbol{X}_0 \in \mathbb{C}^{M \times N}$ where $\boldsymbol{X}$ is arbitrary and $\boldsymbol{X}_0$ has no zero singular values. Substituting $\boldsymbol{Y} = \boldsymbol{X}^*\boldsymbol{X} + \epsilon\boldsymbol{I} \in H_{++}^N$, $\boldsymbol{Y}_0 = \boldsymbol{X}_0^*\boldsymbol{X}_0 + \epsilon\boldsymbol{I} \in H_{++}^N$, and $q = p/2$ into (7), it follows that the function

$$\begin{aligned} &g_p(\boldsymbol{X}; \boldsymbol{X}_0) \\ &= Tr[(\boldsymbol{X}_0^*\boldsymbol{X}_0 + \epsilon\boldsymbol{I})^{\frac{p}{2}} + \frac{p}{2}(\boldsymbol{X}_0^*\boldsymbol{X}_0 + \epsilon\boldsymbol{I})^{\frac{p}{2}-1}(\boldsymbol{X}^*\boldsymbol{X} - \boldsymbol{X}_0^*\boldsymbol{X}_0)] \\ &= \frac{p}{2}Tr[(\boldsymbol{X}_0^*\boldsymbol{X}_0 + \epsilon\boldsymbol{I})^{\frac{p}{2}-1}\boldsymbol{X}^*\boldsymbol{X}] + C(\boldsymbol{X}_0) && (10) \\ &= \frac{p}{2}\left\|\boldsymbol{X}\sqrt{(\boldsymbol{X}_0^*\boldsymbol{X}_0 + \epsilon\boldsymbol{I})^{\frac{p}{2}-1}}\right\|_F^2 + C(\boldsymbol{X}_0), && (11) \end{aligned}$$

where $C(\boldsymbol{X}_0)$ is a term depending only on $\boldsymbol{X}_0$, satisfies the majorization relations

$$g_p(\boldsymbol{X}; \boldsymbol{X}_0) \geq \|\boldsymbol{X}\|_{p,\epsilon}^p, \quad \text{for all} \quad \boldsymbol{X} \in \mathbb{C}^{M \times N} \tag{12}$$

$$g_p(\boldsymbol{X}_0; \boldsymbol{X}_0) = \|\boldsymbol{X}_0\|_{p,\epsilon}^p \tag{13}$$

for all $p \in (0, 1]$. A similar argument using inequality (8) shows that (12) and (13) can be extended to hold for the $p = 0$ case, with the majorizer

$$g_0(\boldsymbol{X}; \boldsymbol{X}_0) = \frac{1}{2}\left\|\boldsymbol{X}\sqrt{(\boldsymbol{X}_0^*\boldsymbol{X}_0 + \epsilon\boldsymbol{I})^{-1}}\right\|_F^2 + C(\boldsymbol{X}_0). \tag{14}$$

Finally, substituting $\boldsymbol{X} = \mathcal{T}(\boldsymbol{x})$ gives the following MM scheme for minimizing the smoothed Schatten-$p$ penalty:

$$\boldsymbol{x}^{(n)} = \arg\min_{\boldsymbol{x}} \|\boldsymbol{A}\boldsymbol{x} - \boldsymbol{b}\|_2^2 + \lambda C_p \|\mathcal{T}(\boldsymbol{x})\sqrt{\boldsymbol{H}_n}\|_F^2 \tag{15}$$

where $C_p = \frac{p}{2}$ if $0 < p \leq 1$, $C_0 = \frac{1}{2}$, and $\boldsymbol{H}_n$ is determined by the previous iterate $\boldsymbol{x}^{(n-1)}$ according to

$$\boldsymbol{H}_n = [\mathcal{T}(\boldsymbol{x}^{(n-1)})^*\mathcal{T}(\boldsymbol{x}^{(n-1)}) + \epsilon\boldsymbol{I}]^{\frac{p}{2}-1}. \tag{16}$$

The above shows that the IRLS-$p$ algorithm (Alg. 1) is a MM algorithm for a smoothed Schatten-p norm penalty with fixed $\epsilon$. To avoid convergence to local minima, we use an iteration dependent $\epsilon_n > 0$ that exponentially decreases to pre-determined minimum value $\epsilon_{\text{min}} > 0$, such that $\epsilon_n = \epsilon_{\text{min}}$ for all iterations $n \geq n'$ for some finite $n'$. Therefore, IRLS-$p$ algorithm presented in this work only behaves as an MM scheme for iterations $n' > n$ where the smoothing parameter is fixed. However, this still guarantees that in the long-run the algorithm converges to a stationary point of (5) with $\epsilon = \epsilon_{\text{min}}$.

## II. Existing algorithms for structured low-rank matrix recovery

### A. Alternating projections algorithms

Given an estimate of the underlying rank $r$ of lifted matrix $\mathcal{T}(\boldsymbol{x}_0)$, one approach to recover $\boldsymbol{x}_0$ from its linear measurements $\boldsymbol{b} = \boldsymbol{A}\boldsymbol{x}_0$ is to solve:

$$\min_{\boldsymbol{x}} \|\boldsymbol{A}\boldsymbol{x} - \boldsymbol{b}\|^2 \quad \text{subject to} \quad \begin{cases} \boldsymbol{X} = \mathcal{T}(\boldsymbol{x}) \\ \operatorname{rank} \boldsymbol{X} \leq r. \end{cases} \tag{17}$$

The alternating projections (AP) algorithm, also known as Cadzow's method after [4], seeks to find a minimum of (17) by alternately projecting onto: (1) the set of matrices with rank less than or equal to $r$, (2) the space of linear structured matrices specified by the range of $\mathcal{T}$, and (3) the data fidelity constraint set. These projections are computed by (1) a rank $r$ truncated SVD (2) an averaging operation determined by the pseudo-inverse $\mathcal{T}^\dagger = (\mathcal{T}^*\mathcal{T})^{-1}\mathcal{T}$, and (3) a least-squares problem, which often has a closed-form solution. Pseudo-code for this approach is shown in Algorithm 1.

**Algorithm 1:** AP algorithm for structured low-rank matrix recovery

Initialize $\boldsymbol{x} \in \mathbb{C}^n$;
**while** *not converged* **do**
  $\boldsymbol{X} = \boldsymbol{U}_r \boldsymbol{\Sigma}_r \boldsymbol{V}_r^*$ where $(\boldsymbol{U}_r, \boldsymbol{\Sigma}_r, \boldsymbol{V}_r) = \text{svd}_r(\mathcal{T}(\boldsymbol{x}))$;
  $\boldsymbol{y} = \mathcal{T}^\dagger(\boldsymbol{X})$;
  $\boldsymbol{x} = (\boldsymbol{I} - \boldsymbol{A}^\dagger \boldsymbol{A})\boldsymbol{y} + \boldsymbol{A}^\dagger \boldsymbol{b}$;
**end**

A novel adaptation of the alternating projections algorithm was proposed in the LORAKS framework [5], which relaxes the rank constraint in (17) by introducing the following (non-convex) functional

$$J_r(\boldsymbol{X}) := \min_{\boldsymbol{T}:\ \operatorname{rank} \boldsymbol{T} \leq r} \|\boldsymbol{X} - \boldsymbol{T}\|_F^2$$

i.e., $J_r(\boldsymbol{X})$ measures the distance to the best rank $r$ approximation of a matrix $\boldsymbol{X}$. This can also be interpreted as a proximal smoothing [6] of the indicator function $I$ given by $I(\boldsymbol{X}) = 0$ if $\operatorname{rank} \boldsymbol{X} \leq r$, and $I(\boldsymbol{X}) = \infty$ otherwise. As a surrogate for (17), [5] proposed to minimize:

$$\min_{\boldsymbol{x}} \|\boldsymbol{A}\boldsymbol{x} - \boldsymbol{b}\|^2 + \lambda J_r(\boldsymbol{X}) \quad \text{subject to} \quad \boldsymbol{X} = \mathcal{T}(\boldsymbol{x}) \tag{18}$$

where $\lambda$ is a regularization parameter. Observe that the objective in (18) approaches the Cadzow formulation (17) as $\lambda \to \infty$. The authors in [5] propose an alternating minimization scheme for solving (18); the pseudo-code for this algorithm is shown in (2), which we call AP-PROX.

**Algorithm 2:** AP-PROX algorithm for structured low-rank matrix recovery [5], [7]

Initialize $\boldsymbol{x} \in \mathbb{C}^n$;
**while** *not converged* **do**
  $\boldsymbol{X} = \boldsymbol{U}_r \boldsymbol{\Sigma}_r \boldsymbol{V}_r^*$ where $(\boldsymbol{U}_r, \boldsymbol{\Sigma}_r, \boldsymbol{V}_r) = \text{svd}_r(\mathcal{T}(\boldsymbol{x}))$;
  $\boldsymbol{x} = (\boldsymbol{A}^*\boldsymbol{A} + \lambda \mathcal{T}^*\mathcal{T})^{-1}(\boldsymbol{A}^*\boldsymbol{b} + \lambda \mathcal{T}^*(\boldsymbol{X}))$;
**end**

### B. Singular value thresholding algorithms

Another common approach to structured low-rank matrix recovery is replace the rank functional with its convex relaxation, the nuclear norm $\|\cdot\|_*$, and solve

$$\min_{\boldsymbol{x}} \|\boldsymbol{A}\boldsymbol{x}-\boldsymbol{b}\|^2 + \lambda\|\boldsymbol{X}\|_* \text{ subject to } \boldsymbol{X} = \mathcal{T}(\boldsymbol{x}), \tag{19}$$

where $\lambda$ is a regularization parameter. Applying the ADMM algorithm [8] to (19) results in the singular value thresholding (SVT) algorithm, originally proposed by [9] in the setting of low-rank matrix completion. See Algorithm 3 for the pseudo-code for this approach. The main drawback of the SVT approach is that it requires a full SVD of the a matrix having dimensions the size of the matrix lifting at every iteration. A well-known workaround to this problem is to exploit the variational characterization of the nuclear norm [10]:

$$\|\boldsymbol{X}\|_* = \min_{\substack{\boldsymbol{U}\in\mathbb{C}^{M\times R}\\ \boldsymbol{V}\in\mathbb{C}^{N\times R}\\ \boldsymbol{X}=\boldsymbol{U}\boldsymbol{V}^*}} \frac{1}{2}\left(\|\boldsymbol{U}\|_F^2 + \|\boldsymbol{V}\|_F^2\right), \tag{20}$$

which holds true provided the inner dimension $R$ of $\boldsymbol{U}$ and $\boldsymbol{V}^*$ satisfies $R \geq \operatorname{rank}\boldsymbol{X}$. Substituting (20)

---

**Algorithm 3:** SVT algorithm for structured low-rank matrix recovery

---

Initialize $\boldsymbol{x}\in\mathbb{C}^n$, $\boldsymbol{L}=\boldsymbol{0}\in\mathbb{C}^{M\times N}$, and choose $\mu>0$;
**while** *not converged* **do**
  $\boldsymbol{X} = \boldsymbol{U}\max(\boldsymbol{\Sigma}-\frac{1}{\mu}\boldsymbol{I},0)\boldsymbol{V}^*$ where $(\boldsymbol{U},\boldsymbol{\Sigma},\boldsymbol{V}) = \operatorname{svd}(\mathcal{T}(\boldsymbol{x})+\boldsymbol{L})$;
  $\boldsymbol{x} = (\mu\mathcal{T}^*\mathcal{T}+\lambda\boldsymbol{A}^*\boldsymbol{A})^{-1}[\mu\mathcal{T}^*(\boldsymbol{X}-\boldsymbol{L})+\lambda A^*\boldsymbol{b}]$;
  $\boldsymbol{L} = \boldsymbol{L}+\mathcal{T}(\boldsymbol{x})-\boldsymbol{Z}$
**end**

---

into (19) gives the following constrained optimization problem:

$$\min_{\boldsymbol{x},\boldsymbol{U},\boldsymbol{V}} \|\boldsymbol{A}\boldsymbol{x}-\boldsymbol{b}\|^2 + \frac{\lambda}{2}\left(\|\boldsymbol{U}\|_F^2+\|\boldsymbol{V}\|_F^2\right)$$
$$\text{subject to } \boldsymbol{U}\boldsymbol{V}^* = \mathcal{T}(\boldsymbol{x}), \tag{21}$$

which can similarly be solved using ADMM. The resulting algorithm has the same structure as the original SVT approach, except the singular value thresholding step is replaced with matrix inversion steps. We call this approach the SVT+UV algorithm; see Algorithm 4 for the pseudo-code. While the SVT+UV algorithm

---

**Algorithm 4:** SVT+UV algorithm for structured low-rank matrix recovery [11]–[13]

---

Initialize $\boldsymbol{x}\in\mathbb{C}^n$, $\boldsymbol{L}=\boldsymbol{0}\in\mathbb{C}^{M\times N}$, $\boldsymbol{V}\neq\boldsymbol{0}\in\mathbb{C}^{N\times R}$, and choose $\mu>0$;
**while** *not converged* **do**
  $\boldsymbol{C} = \mu(\mathcal{T}(\boldsymbol{x})+\boldsymbol{L})$;
  $\boldsymbol{U} = (\boldsymbol{C}\boldsymbol{V})/(\boldsymbol{I}_R+\mu\boldsymbol{V}^*\boldsymbol{V})$;
  $\boldsymbol{V} = (\boldsymbol{C}^*\boldsymbol{U})/(\boldsymbol{I}_R+\mu\boldsymbol{U}^*\boldsymbol{U})$;
  $\boldsymbol{X} = \boldsymbol{U}\boldsymbol{V}^*$;
  $\boldsymbol{x} = (\mu\mathcal{T}^*\mathcal{T}+\lambda\boldsymbol{A}^*\boldsymbol{A})^{-1}[\mu\mathcal{T}^*(\boldsymbol{X}-\boldsymbol{L})+\lambda\boldsymbol{A}^*\boldsymbol{b})]$;
  $\boldsymbol{L} = \boldsymbol{L}+\mathcal{T}(\boldsymbol{x})-\boldsymbol{X}$;
**end**

---

is more computationally efficient than SVT, it can have several practical drawbacks. First, by introducing the variables $\boldsymbol{U},\boldsymbol{V}$, the objective in (II-B) becomes non-convex due to the constraint $\mathcal{T}(\boldsymbol{x}) = \boldsymbol{U}\boldsymbol{V}^*$. An essentially similar problem has been considered in [11] and, despite non-convexity, the algorithm has been observed to converge provided the rank $r$ of the global minimizer is sufficiently smaller than rank parameter $R$ [11]. Related algorithms have been proven to converge to the global minimum under this condition as well [10], [14], and it is likely this analysis can be carried over to the SVT+UV algorithm. However, the algorithm is less stable when $R < r$, and can require special initialization for good convergence in this case. Finally, the algorithm still has significant memory demands, since it requires storing a dual variable having dimensions of the lifted matrix ($\boldsymbol{L}$ in Algorithm 4).